\documentclass{iopart}
\usepackage{iopams}
\usepackage{amsthm}
\usepackage{times}
\usepackage{url}
\usepackage[matrix,arrow]{xy}
\usepackage{graphicx}

\newcommand{\bfi}{\bfseries\itshape}
\def\eqref#1{(\ref{#1})}  
\def\operatorname#1{\mathrm{#1}}
\def\intprod{\mathbin{\hbox to 6pt{ \vrule height0.4pt width5pt depth0pt
\kern-.4pt \vrule height6pt width0.4pt depth0pt\hss}}}

\newtheorem{theorem}{Theorem}[section]

\newtheorem{lemma}[theorem]{Lemma}

\begin{document}

\title[Discrete Routh Reduction]{Discrete Routh Reduction}

\author{Sameer M. Jalnapurkar$^1$, Melvin Leok$^2$, Jerrold E. Marsden$^3$ and Matthew West$^4$}

\address{$^1$ Department of Mathematics, Indian Institute of Science, Bangalore, India.}
\address{$^2$ Department of Mathematics, University of Michigan, East Hall, 530 Church Street, Ann Arbor, MI 48109-1043, USA.}
\address{$^3$ Control and Dynamical Systems 107-81, California Institute of Technology, Pasadena, CA 91125-8100, USA.}
\address{$^4$ Department of Aeronautics and Astronautics, Stanford
University, Stanford, CA 94305-4035, USA.}
\ead{mleok@umich.edu}
\begin{abstract}
This paper develops the theory of abelian Routh reduction  for discrete mechanical systems and applies it to the variational integration of mechanical systems with abelian symmetry. The reduction of variational Runge-Kutta discretizations is considered, as well as the extent to which symmetry reduction and discretization commute. These reduced methods allow the direct simulation of dynamical features such as relative equilibria and relative periodic orbits that can be obscured or difficult to identify in the unreduced dynamics.

The methods are demonstrated  for the dynamics of an Earth orbiting satellite with a non-spherical $J_2$ correction, as well as the double spherical pendulum. The $J_2$ problem is interesting because in the unreduced picture, geometric phases inherent in the model and those due to numerical discretization can be hard to distinguish, but this issue does not appear in the reduced algorithm, where one can directly observe interesting dynamical structures in the reduced phase space (the cotangent bundle of shape space), in which  the geometric phases have been removed.

The main feature of the double spherical pendulum example is that it has a nontrivial magnetic term in its reduced symplectic form. Our method is still efficient as it can directly handle the essential non-canonical nature of the symplectic structure. In contrast, a traditional symplectic method for canonical systems could require repeated coordinate changes if one is evoking Darboux' theorem to transform the symplectic structure into canonical form, thereby incurring additional computational cost. Our method allows one to design reduced symplectic integrators in a natural way, despite the noncanonical nature of the symplectic structure.
\end{abstract}

\section{Introduction} \label{s-intr}
This paper addresses reduction theory for discrete mechanical systems with abelian symmetry groups and its relation to variational integration. To establish the setting of the problem, a few aspects of the continuous theory are first recalled (see \cite{MaRa1999} for general background).
\medskip 

\noindent {\bf Continuous Reduction Theory.} Consider a mechanical system with configuration manifold $Q$ and a symmetry group $G$ (with Lie algebra $\mathfrak{g}$) acting freely and properly on $Q$ and hence, by cotangent lift on $T^{\ast}Q$, with corresponding (standard, equivariant) momentum map $\mathbf{J}: T^{\ast}Q \rightarrow \mathfrak{g}^\ast$. Recall from reduction theory (see \cite{MaWe1974, Marsden1992}, and references therein) that, under appropriate regularity and nonsingularity conditions, the flow of a $G$-invariant Hamiltonian $H: T^{\ast}Q \rightarrow \mathbb{R}$ naturally induces a Hamiltonian flow on the reduced space $P _\mu = \mathbf{J}^{-1}(\mu)/ G _\mu$, where $G _\mu$ is the isotropy subgroup of a chosen point $\mu \in \mathfrak{g}^\ast$. In the abelian case, if one chooses a connection $\mathfrak{A}$ on the principal bundle $Q \rightarrow Q/G$, then $P_\mu$ is symplectically isomorphic to $T^{\ast}(Q/G)$ carrying the canonical symplectic structure modified by magnetic terms; that is terms induced from the $\mu$-component of the curvature of $\mathfrak{A}$.

The Lagrangian version of this theory is also well-developed. In the abelian case, it goes by the name of  {\it Routh reduction}, (see, for instance \cite{MaRa1999}, \S8.9).  The reduced equations are again equations on $T(Q/G)$ and are obtained by dropping the variational principle, expressed in terms of the Routhian, from $Q$ to $Q/G$. The nonabelian version of this theory was originally developed in \cite{MaSc1993a, MaSc1993b}, with important contributions and improvements given in \cite{JaMa2000, MaRaSc2000}.

Of course, reduction has been enormously important for many topics in mechanics, such as stability and bifurcation of relative equilibria, integrable systems, etc. We need not review the importance of this process here as it is extensively documented in the literature.

\medskip 

\noindent {\bf Purpose, Main Results, and Examples.} This paper presents the theory and illustrative numerical implementation for the reduction of discrete mechanical systems with abelian symmetry groups. The discrete reduced space has a similar structure as in the continuous theory, but the curvature will be taken in a discrete sense. The paper studies two examples in detail, namely satellite dynamics in the presence of the bulge of the Earth (the $J _2$ effect) and the double spherical pendulum (which has a nontrivial magnetic term). In each case the benefit of studying the numerics of the reduced problem is shown. Roughly, the reduced computations reveal dynamical structures that are hard to pick out in the unreduced dynamics in a way that is reminiscent of the phenomena of pattern evocation, as in \cite{MaSc1995, MaScWe1996}. Another interesting application of the theory is that of orbiting multibody systems, studied in \cite{SaShMc2003,SaBlMc2004}.

We refer to \cite{MaWe2001} for a review of discrete mechanics, its numerical implementation, some history, as well as references to the literature. The value of geometric integrators has been documented in a number of references, such as~\cite{HaLuWa2001}. In the present paper, we shall focus, to be specific, on discrete Euler--Lagrange and variational symplectic Runge-Kutta schemes and their reductions. One could, of course, use other schemes as well, such as Newmark, St\"ormer-Verlet, or Shake schemes. However, we wish to emphasize that {\it without theoretical guidelines, coding algorithms for the reduced dynamics need not be a routine procedure since the reduced equations are not in canonical form because of nontrivial magnetic terms}. For example, using Darboux' theorem to put the structure into canonical form so that standard algorithms can be used is not practical. We also remind the reader that there are real advantages to taking the variational approach to the construction of symplectic integrators. For example, as in \cite{LeMaOrWe2003}, the variational approach provides the design flexibility to take different time steps at spatially different points in an asynchronous way and still retain all the advantages of symplecticity even though the algorithms are not strictly symplectic in the naive sense; such an approach is well-known to be useful in molecular systems, for instance.

\medskip 

\noindent {\bf Motivation for Discrete Reduction.} Besides its considerable theoretical interest, there are several practical reasons for carrying out discrete Routh reduction. These are as follows:
\noindent\begin{enumerate}
\item Features that are clear in the reduced dynamics, such as relative equilibria and relative periodic orbits, can be obscured in the unreduced dynamics, and appear more complicated through the process of reconstruction and associated geometric phases. This is related to the phenomenon of {\it pattern evocation} that is an important practical feature of many examples, such as the double spherical pendulum \cite{MaSc1995, MaScWe1996} and the stepping pendulum \cite{HoLy2002}. Going to a suitable (but non-obvious) rotating frame can ``evoke'' such phenomena (see the movie at \url{http://www.cds.caltech.edu/~marsden/research/demos/movies/Wendlandt/pattern.mpg}). This is essentially a window to the reduced dynamics, which the theory in the present paper allows one to compute directly.

\item While directly studying the reduced dynamics can yield some benefits, it can be difficult to code using traditional methods. In particular, the presence of magnetic terms in the reduced symplectic form, as is the case with the double spherical pendulum, means that traditional symplectic methods for canonical systems do not directly apply; if one attempts to do so, it may result (and has in the literature) in many inefficient coordinate changes when evoking Darboux' theorem to put things into canonical form.

\item Although simulating the reduced dynamics involves an initial investment of time in computing geometric quantities symbolically, these additional terms do not appreciably affect the sparsity of the system of equations to be solved. As such, direct coding of the reduced algorithms can be quite efficient, due to its reduced dimensionality. \end{enumerate}

\medskip 

\noindent {\bf Two Obvious Generalizations.} The free and proper assumption that we make on the group action means that we are dealing with {\it nonsingular}, that is, {\it regular} reduction (see \cite{OrRa2004} for the general theory of singular reduction and references to the literature). It would be interesting to extend the work here to the case of singular reduction but already the regular case is nontrivial and interesting. While our examples have singular points and the dynamics near these points is interesting, there is no attempt to study this aspect in the present paper.

Secondly, it would be interesting to generalize the present work to the case of nonabelian groups and to develop a discrete version of nonabelian Routh reduction, (as in \cite{JaMa2000, MaRaSc2000}). We believe that such a generalization will require the further development of the theory of discrete connections, which is currently part of the research effort on discrete differential geometry (see \cite{Leok2004}, and references therein.) Other future directions are discussed in the conclusions.

\medskip 

\noindent {\bf Other Discrete Reduction Results.} We briefly summarize some related results that have been obtained in the area of reduction for discrete mechanics. First of all, there is the important case of {\it discrete Euler-Poincar\'e and Lie-Poisson reduction} that were obtained in \cite{BoLoSu1998, MaPeSh1999,MaPeSh2000}. This theory is appropriate for rigid body mechanics, for instance.

Another important case is that of {\it discrete semidirect product reduction} that was obtained in \cite{BoSu1999, BoSu1999a} and applied to the case of the heavy top, with interesting links to discrete elastica. This case is of interest in the present study since the heavy top, as with the general theory of semidirect product reduction (see \cite{MaRaWe1984, HoMaRa1998}), one can view the $S ^1$ reduction of this problem as Routh reduction. Linking these two approaches is an interesting topic for future research.
\medskip 

\noindent {\bf Outline.}
After recalling the notation from continuous reduction theory, \S\ref{s-disRed} develops discrete reduction theory, derives a reduced variational principle, and proves the symplecticity of the reduced flow. The relationship between continuous- and discrete-time reduction is also discussed. How the variational (and hence symplectic) Runge-Kutta algorithm induces a reduced algorithm in a natural way is shown in \S\ref{s-srkRed}. In \S\ref{s-summary} we put together in a coherent way the main theoretical results of the paper up to that point. In \S\ref{s-examples} the numerical example of satellite dynamics about an oblate Earth is given, and in \S\ref{s-examples-dsp}, the example of the double spherical pendulum, which has a non-trivial magnetic term, is given. Lastly, in \S\ref{s-computational}, we address  some computational and efficiency issues.

\section{Discrete Reduction}\label{s-disRed} In this section, it is assumed that the reader is familiar with continuous reduction theory as well as the theory of discrete mechanics; reference is made to the relevant parts of the literature as needed. It will be useful to first recall some facts about discrete mechanical systems with symmetry (see, for instance, \cite{MaWe2001} for proofs.)

\medskip 

\noindent {\bf Discrete Mechanical Systems with Symmetry.} Let $G$ be a Lie group (which shortly will be assumed to be abelian) that acts freely and properly (on the left) on a configuration manifold $Q$. Given a discrete Lagrangian $L_d: Q \times Q \rightarrow \mathbb{R}$  that is invariant under the diagonal action of $G$ on $Q\times Q$, the corresponding {\bfi discrete momentum map} $\mathbf{J}_d :Q\times Q \to \mathfrak{g}^*$ is defined by
\begin{equation}
\label{discrete-mom-map}
\mathbf{J}_d (q_0,q_1)\cdot \xi = D_2L_d(q_0,q_1)\cdot \xi_Q(q_1),
\end{equation}
where $D_2$ denotes the derivative in the second slot and where $\xi_Q$ is the infinitesmal generator associated with $\xi \in \mathfrak{g}$. The map $\mathbf{J}_d$ is equivariant with respect to the diagonal action of $G$ on $Q \times Q $ and the coadjoint action on $\mathfrak{g}^\ast$. The {\bfi discrete Noether theorem} states that the discrete momentum is conserved along solutions of the DEL (Discrete Euler--Lagrange) equations,
\begin{equation}
\label{DEL-equations}
     D_2L_d(q_{k-1},q_k) + D_1L_d(q_k,q_{k+1}) = 0.
\end{equation}
 Notice that
\[
\mathbf{J}_d (q_0,q_1)\cdot\xi = \mathbf{J}(D_2 L_d (q_0,q_1))\cdot \xi ,
\]
where $\mathbf{J}:T^{\ast}Q\to \mathfrak{g}^*$ is the momentum map on $T^{\ast}Q$; i.e.,
$
\mathbf{J}_d =\mathbf{J}\circ \mathbb{F}L_d,
$
where $\mathbb{F}L_d = D_2 L_d : Q\times Q \to T^{\ast}Q$ is the discrete Legendre transform. Thus, for $\mu \in \mathfrak{g}^\ast$, we have $\mathbb{F}L_d\left(\mathbf{J}_d ^{-1}(\mu)\right) \subset \mathbf{J}^{-1}(\mu)$. The symplectic algorithm (usually called the {\it position-momentum form of the algorithm}) obtained on $T^{\ast}Q$ from that on $Q \times Q $ via the discrete Legendre transform thus preserves the standard momentum map $\mathbf{J}$.

There will be a standing assumption in this paper, namely that the given discrete Lagrangian $L_d$ is {\bfi regular}; that is, for a point $(q,q ) \in Q \times Q$ on the diagonal, the iterated derivative $D_2 D_1 L_d(q,q): T_{q} Q \times T_{q}Q \rightarrow \mathbb{R}$ is a non-degenerate bilinear form. By the implicit function theorem, this implies that a point $(q_{k-1},q_k)$ near the diagonal and the DEL equations \eqref{DEL-equations} uniquely determines the subsequent point $q _{k + 1}$ in a neighborhood of the diagonal in $Q \times Q $ (or, if one prefers, for small time steps); in other words, the DEL algorithm is well defined by the  DEL equations. Regularity also implies that the discrete Legendre transformation $\mathbb{F}L_d = D_2 L_d : Q \times Q \to T^{\ast}Q$ is a local diffeomorphism from a neighborhood of the diagonal in $Q \times Q $ to a neighborhood in $T^{\ast}Q$. For a detailed discussion, see \cite{MaWe2001}.
\medskip 

\noindent {\bf Reconstruction.} The following lemma gives a basic result on the reconstruction of discrete curves in the configuration  manifold $Q$ from those in {\bfi shape space}, defined to be $S =Q/G$. The Lemma is similar to its continuous counterpart, as in, for example, \cite{JaMa2000}, Lemma 2.2. The natural projection to the quotient will be denoted $\pi_{Q,G}: Q \rightarrow Q/G $; $q \mapsto x = [q] _G$ (the equivalence class of $q \in Q $). Let $\operatorname{Ver}(q)$ denote the {\bfi vertical space} at $q$; namely the space of all vectors at the point $q$ that are infinitesimal generators $\xi _Q (q) \in T _q Q $; or in other words, the tangent space to the group orbit through $q$. We say that the discrete Lagrangian $L_d$ is {\bfi group-regular} if the bilinear map $D_2 D_1 L_d(q,q): T_{q} Q \times T_{q}Q \rightarrow \mathbb{R}$ restricted to the subspace $\operatorname{Ver}(q) \times \operatorname{Ver}(q)$ is nondegenerate. {\it In addition to regularity, we shall make group-regularity a standing assumption in the paper.}
The following result is fundamental for what follows.
\begin{lemma}[Reconstruction Lemma]
\label{reconstruct-lemma}
Fix $\mu \in \mathfrak{g}^\ast$ and let $ x_0, x _1, \ldots,x_n$ be a sufficiently closely spaced discrete curve in $S$. Let $q_0, q_1  \in Q $ be such that $[q_0] _G = x_0, [q_1] _G = x_1$ and $\mathbf{J}_d (q_0, q_1) = \mu $. Then there is a unique closely spaced discrete curve $ q_1, q_2,   \ldots, q_n$ such that  $[q_k] _G = x_k$ and $\mathbf{J}_d (q_{k - 1},q_{k})=\mu$, for $k=1, 2,\ldots,n$.
\end{lemma}

\begin{proof} We must construct a point $q _2$ close to $q _1$ such that $[q_2] _G = x_2$ and $\mathbf{J}_d (q_1,q_{2})=\mu$; the construction of the subsequent points $q _3, \ldots q _n$ then proceeds in a similar fashion.

To do this, pick a local trivialization of the bundle $\pi_{Q,G}: Q \rightarrow Q/G$ where locally $Q = S \times G $ and write points in this trivialization as  $q_0 = (x_0, g_0)$ and $q_1 = (x_1, g_1)$, etc. Given the points $q_0 = (x_0, g_0), q_1 = (x_1, g_1)$ with
$\mathbf{J}_d(q_0, q_1) = \mu$, we seek a near identity group element $k \in G$ such that $q_2 : = (x_2, k g_1) $ satisfies $\mathbf{J}_d (q_{1},q_{2}) = \mu$. According to equation \eqref{discrete-mom-map}, this means that we must satisfy the 
condition
$D_2L_d(q_1,q_2) \cdot  \xi_Q(q_2) = \left\langle \mu, \xi \right\rangle,
$  
for all $\xi \in \mathfrak{g}$. In the local trivialization, this reads
 \begin{equation} \label{reconstruction_equation}
 D_2L_d((x_1, g_1),(x_2, k g_1)) \cdot ( 0,  TR _{k g _1} \xi) = \left\langle \mu, \xi \right\rangle,
 \end{equation}
where $R _g$ is right translation on $G$ by $g$. Consider solving  the equation 
 \begin{equation} \label{reconstruction_eqn2}
 D_2L_d((\bar{x}_1, \bar{g}_1),(\bar{x}_2, k \bar{g}_1)) \cdot ( 0,  TR _{k \bar{g} _1} \xi) = \left\langle \mu, \xi \right\rangle,
 \end{equation} 
for $k$ as a function of the variables $ \bar{g}_1, \bar{x}_1, \bar{x}_2$ with $\mu$ fixed. By assumption, there is a  solution for the case $\bar{x}_1= x_0$, $\bar{x} _2 = x _1$ and $\bar{g} _1 = g _0$, namely $k = k_0 = g _1 g _0 ^{-1}$ (a near identity group element). The implicit function theorem shows that when the point $g_0, x_0, x_1$ is replaced by the nearby point $g_1, x_1, x_2$, there will be a unique solution for $k$ near $k _0$ provided that the derivative of the defining relation \eqref{reconstruction_equation} with respect to $k$ at the identity is invertible, which is true by  group-regularity. Since group regularity is $G$-invariant, the above argument remains valid as $k_i$ drifts from the identity.
\end{proof}

Note that the above Lemma makes no hypotheses about the sequences $x _n$ or $q _n$ satisfying any discrete evolution equations.

To carry out reconstruction in the continuous case, in addition to the requirements that the lifted curve in $TQ$ lie on the momentum surface, and that it projects to the reduced curve $x(t)\in S$ under $\pi_{Q,G}$, one also requires that it be second-order, which is to say that it is of the form $(q(t), \dot{q}(t))$. If a connection is given, then the lifted curve is obtained by integrating the {\it reconstruction equation}--again, see \cite{JaMa2000} for details. The discrete analogue of the second-order curve condition is explained as follows. Consider a given discrete curve as a sequence of points, $(x_0,x_1),(x_1,x_2),\ldots,(x_{n-1},x_n)$ in $S\times S$. Lift each of the points in $S\times S$ to the momentum surface $\mathbf{J}_d^{-1}(\mu)\subset Q\times Q$. This yields the sequence, $(q^0_0,q^0_1),(q^1_0,q^1_1),\ldots,(q^{n-1}_0,q^{n-1}_1)$, which is unique up to an overall diagonal group action. The discrete analogue of the second-order curve condition is that this sequence in $Q\times Q$ defines a discrete curve in $Q$, which corresponds to requiring that $q^k_1=q^{k+1}_0$, for $k=0,\ldots,n-1$, which is clearly possible in the context of the reconstruction lemma.

Discrete reconstruction naturally leads to issues of discrete geometric phases, and it would be interesting to express the discrete geometric phase in terms of the discrete curvature on shape space; this will surely involve some ideas from the currently evolving subject of discrete differential geometry and so we do not attempt to push this idea further at this point.

While many of the computations we present in this paper are in the setting of local trivializations, the results are valid  globally through the construction given below.
\medskip 

\noindent {\bf Identification of the Quotient Space.} Now assume that $G $ is abelian so that $G = G _\mu$ acts on $\mathbf{J}_d^{-1}(\mu)$ and so that the quotient space $\mathbf{J}_d^{-1}(\mu)/G $ makes sense. We assume the above regularity hypotheses and  freeness and properness of the action of $G$ so that this quotient is a smooth manifold. It is clear that the map $\varphi _\mu: \mathbf{J}_d^{-1}(\mu)/G  \rightarrow S \times S$ given by $[ (q, q^{\prime})] _G \rightarrow (x, x^{\prime}) $ is well-defined, where the square brackets denotes the equivalence class with respect to the given $G$ action, and where $x = [q] _G$.  The argument given in the reconstruction lemma shows that for a point $(q_0, q_1) \in \mathbf{J}_d^{-1}(\mu)$,  $\varphi_\mu$ is a local diffeomorphism in a neighborhood of the point $[(q_0, q_1)] _G$. In fact, the uniqueness part of that argument shows that for two nearby points $(q_1, q_2)$ and $(q^{\prime}_1, q^{\prime}_2)$ in $\mathbf{J}_d^{-1}(\mu)$, if  $q_1 = g_1 q^{\prime}_1$ and $q_2 = g_2 q^{\prime}_2$, then $g_1 = g_2$. Thus, there is a neighborhood $U$ of a given a chain of closely spaced points lying in $\mathbf{J}_d^{-1}(\mu)$ with this property. Saturating this neighborhood with the group action, we can assume that $U$ is $G$-invariant. Restricted to $U$, $\varphi _\mu$ becomes a diffeomorphism to a neighborhood of the diagonal of $S \times S$.  

Assume, as above, that $L_d : Q \times Q \to \mathbb{R}$ is a discrete Lagrangian that is invariant under the action of an abelian Lie group $G$ on $Q \times Q$. In view of the preceding discussion, $L _d$ restricted to $\mathbf{J}_d^{-1}(\mu)$ (and in the neighborhood of a given chain of points in this set) induces a well defined function $ \hat{L }_d ( x _0, x _1)$ of pairs of points $( x _0, x _1)$ in $S \times S $. This {\bfi discrete reduced Lagrangian} will play an important role in what follows.

\medskip 

\noindent {\bf Discrete Reduction.}
Let $\mathbf{q}:=\{q_0,\ldots,q_n\}$ be a solution of the
discrete Euler-Lagrange (DEL) equations. Let the value of the discrete momentum along this trajectory be $\mu$. Let $x_i=[q_i] _G$, so that $\mathbf{x}:=\{x_0,\ldots,x_n\}$ is a discrete shape space trajectory. Since $\mathbf{q}$ satisfies the discrete variational principle, it is appropriate to ask if there is a reduced variational principle satisfied by $\mathbf{x}$.

An important issue in dropping the discrete variational principle to the shape space is whether we require that the varied curves are constrained to lie on the level set of the momentum map. The constrained approach is adopted in \cite{JaMa2000}, and the unconstrained approach is used in \cite{MaRaSc2000}. In the rest of this section, we will adopt the unconstrained approach of \cite{MaRaSc2000}, and will show that the variations in the discrete action sum, without assuming the variations at the endpoints vanish and evaluated on a solution of the discrete Euler--Lagrange equations, depends only on the quotient variations, and therefore drops to the shape space without constraints on the variations.

If $\mathbf{q}$ is a solution of the DEL equations, the interior terms in the variation of the discrete action sum vanish, leaving  only the boundary terms; that is,
\begin{equation}
\delta \sum_{k=0}^{n-1} L_d(q_k,q_{k+1})=
D_1 L_d(q_0,q_1)\cdot \delta q_0 +D_2 L_d(q_{n-1},q_n)\cdot \delta q_n.
\label{e-dsLhatA}
\end{equation}

Given a  principal connection   $\mathfrak{A}$ on $Q$, there is a horizontal-vertical split of each tangent space to $Q$ denoted $v _q = \operatorname{hor} \, v _q + \operatorname{ver} \, v _q$ for $v _q \in T_qQ$. Thus,
\begin{equation*}
\fl\qquad D_2 L_d(q_{n-1},q_n)\cdot \delta q_n =
D_2 L_d(q_{n-1},q_n)\cdot \operatorname{hor}\,\delta q_n+
D_2 L_d(q_{n-1},q_n)\cdot \operatorname{ver}\,\delta q_n .
\end{equation*}
As in continuous Routh reduction, we will rewrite the terms involving vertical variations using the fact that we are on a level set of $\mathbf{J}_d$. Namely, write the vertical variation as $\operatorname{ver}\,\delta q_n = \xi_Q(q_n)$, where $\xi= \mathfrak{A}(\delta q_n)$ and use the definition \eqref{discrete-mom-map} of $\mathbf{J}_d$ to give
\begin{equation}
\eqalign{
\fl\qquad D_2 L_d(q_{n-1},q_n)\cdot \operatorname{ver}\,\delta q_n &=
D_2 L_d(q_{n-1},q_n)\cdot \xi_Q (q_n)=\mathbf{J}_d (q_{n-1},q_n)\cdot \xi\\
&= \langle \mu,\xi \rangle= \langle \mu,\mathfrak{A}(\delta q_n)\rangle =\mathfrak{A}_\mu
(q_n)\cdot \delta q_n.}
\label{e-crux}
\end{equation}
Thus, the boundary terms can be expressed as
\begin{eqnarray}
\fl\qquad D_2 L_d(q_{n-1},q_n)\cdot \delta q_n & = &
D_2 L_d(q_{n-1},q_n)\cdot \operatorname{hor}\,\delta q_n +
\mathfrak{A}_\mu (q_n)\cdot \delta q_n ,\\
\fl\qquad D_1 L_d(q_0,q_1)\cdot \delta q_0 & = &
D_1 L_d(q_0,q_1)\cdot \operatorname{hor}\,\delta q_0 -
\mathfrak{A}_\mu (q_0)\cdot \delta q_0 ,
\end{eqnarray}
and so \eqref{e-dsLhatA}, the variation of the discrete action sum, becomes
\begin{eqnarray}
 \hspace{- 0.5in} \delta \sum_{k=0}^{n-1} L_d(q_k,q_{k+1}) & = & 
D_1 L_d(q_0,q_1)\cdot \operatorname{hor}\,\delta q_0 +
D_2 L_d(q_{n-1},q_n)\cdot \operatorname{hor}\,\delta q_n \nonumber \\[-6pt]
& &   +  \mathfrak{A}_\mu (q_n)\cdot \delta q_n - \mathfrak{A}_\mu
(q_0)\cdot \delta q_0 , 
\label{e-dsLhatB}
\end{eqnarray}
when restricted to solutions of the discrete Euler-Lagrange equations.

Motivated by the preceding equation, introduce the 1-form $\mathcal{A}$ on $Q\times Q$ defined by
\begin{equation}\label{e-defCalA}
\mathcal{A}=\pi_2^* \mathfrak{A}_\mu - \pi_1^* \mathfrak{A}_\mu,
\end{equation}
where $\pi_1, \pi_2 : Q\times Q\to Q$ are projections onto the first and the
second components respectively. This allows us to expand the boundary terms involving $\mathfrak{A}_\mu$ into a telescoping sum, and rewrite \eqref{e-dsLhatB} in terms of the $1$-form $\mathcal{A}$ as
\begin{eqnarray}
\hspace{- 0.5in}  \sum_{k=0}^{n-1} (DL_d-\mathcal{A})(q_k,q_{k+1})\cdot
(\delta q_k,\delta q_{k+1}) & =& D_1 L_d(q_0,q_1)\cdot
\operatorname{hor}\,\delta q_0 \nonumber \\[-6pt]
&&  + D_2 L_d(q_{n-1},q_n)\cdot \operatorname{hor}\,\delta q_n .
\label{e-dvp}
\end{eqnarray}
We  now drop \eqref{e-dvp} to the reduced space $S\times S$. Consider the projection maps $\pi:Q\times Q\rightarrow (Q\times Q)/G$, and $\pi_{\mu,d}:\mathbf{J}_d^{-1}(\mu)\rightarrow \mathbf{J}_d^{-1}(\mu)/G$, and the inclusion maps, $\iota_{\mu,d}:\mathbf{J}_d^{-1}(\mu)\hookrightarrow Q\times Q$, and $\tilde{\iota}_{\mu,d}:\mathbf{J}_d^{-1}(\mu)/G\hookrightarrow (Q\times Q)/G$. Then clearly, the following diagram commutes,
\[\xymatrix{
\mathbf{J}_d^{-1}(\mu) \; \ar@{^{(}->}[r]^{\iota_{\mu,d}} \ar[d]_{\pi_{\mu,d}} &
Q\times Q \ar[d]^{\pi}\\
\mathbf{J}_d^{-1}(\mu)/G \; \ar@{^{(}->}[r]^{\tilde{\iota}_{\mu,d}} & (Q\times Q)/G \\
}\]
By the $G$-invariance, $L _d$ drops to a function $\tilde L_d$ on the quotient $(Q\times Q)/G$, so that $L_d=\tilde L_d \circ \pi$. The pullback (in this case, the restriction) of $\tilde{L}_d$ to $\mathbf{J}_d^{-1}(\mu)/G$ is called the {\bfi discrete reduced Lagrangian} and is denoted $\hat{L}_d$. Thus,  $\hat{L}_d = \tilde{L}_d \circ \tilde{\iota}_{\mu,d}$; identifying $\mathbf{J}_d^{-1}(\mu)/G$ with $S \times S $, this definition agrees with $\hat{L}_d$ as defined earlier.

\begin{lemma}\label{l-encapsulate-legality}
The $1$-form $\mathcal{A}$ on $Q\times Q$ restricted to $\mathbf{J}_d^{-1}(\mu)$, drops to a $1$-form $\hat{\mathcal{A}}$ on $\mathbf{J}_d^{-1} (\mu)/ G $ and induces, (for closely spaced points), via the map $\varphi _\mu$, a $1$-form that we denote by the same letter, on  $S\times S$. Similarly, the $1$-form $DL_d - \mathcal{A}$ on $Q\times Q$ restricted to the momentum level set $\mathbf{J}_d^{-1}(\mu)$ then drops to the $1$-form  $D\hat{L}_d-\hat{\mathcal{A}}$ on  on $\mathbf{J}_d^{-1} (\mu)/ G $ and induces, via the map $\varphi _\mu$, a $1$-form that we denote by the same letter, on  $S\times S$. 
\end{lemma}
\begin{proof}
The equivariance of the projections $\pi_i$ with respect to the diagonal action on $Q\times Q$ and the given action on $Q$, together with the invariance of the 1-form $\mathfrak{A}_\mu$ on $Q$ imply that $\mathcal{A}$ is invariant. Since $\mathfrak{A}_\mu$ vanishes on vertical vectors for the bundle $Q\rightarrow Q/G$, it follows that $\mathcal{A}$ vanishes on vertical vectors for the bundle $Q\times Q\rightarrow (Q\times Q)/G$. Therefore, there is a 1-form $\tilde\mathcal{A}$ on $(Q\times Q)/G$ such that $\mathcal{A}=\pi^* \tilde\mathcal{A}$.

Since $L_d=\tilde L_d \circ \pi$, and the exterior derivative commutes with pull-back, it follows that $dL_d=\pi^* d\tilde L_d$. From $\pi\circ \iota_{\mu,d}=\tilde{\iota}_{\mu,d}\circ \pi_{\mu,d}$, we get $\iota_{\mu,d}^* \mathcal{A}=\iota_{\mu,d}^*\pi^* \tilde{A} = \pi_{\mu,d}^* \tilde\iota_{\mu,d}^* \tilde{\mathcal{A}}$. Thus, the 1-form $\mathcal{A}$ restricted to $\mathbf{J}_d ^{-1}(\mu)$ drops to the 1-form $\hat{\mathcal{A}}=\tilde{\iota}_{\mu,d}^* \tilde{\mathcal{A}}$ on $\mathbf{J}_d^{-1}(\mu)/G$.
Similarly, $\iota_{\mu,d}^* dL_d =\pi_{\mu,d}^* \tilde{\iota}_{\mu,d}^* d\tilde{L}_d$ and so $dL_d$ restricted to $\mathbf{J}_d ^{-1}(\mu)$ drops to the 1-form $\tilde{\iota}_{\mu,d}^* d\tilde{L}_d = d \hat{L}_d $ on $\mathbf{J}_d^{-1}(\mu)/G$. 
These $1$-forms push forward under the map $\varphi _\mu : \mathbf{J}_d^{-1}(\mu)/G \rightarrow S \times S $ in the manner that was explained earlier. 
\end{proof}

With the preceding lemma, and equation \eqref{e-dvp}, we conclude that
\begin{eqnarray}
\hspace{-0.5in} \sum_{k=0}^{n-1} (D\hat{L}_d-\hat{\mathcal{A}})(x_k,x_{k+1})\cdot
(\delta x_k,\delta x_{k+1}) & = & D_1 L_d(q_0,q_1)\cdot
\operatorname{hor}\,\delta q_0 \nonumber \\[-6pt]
&& + D_2 L_d(q_{n-1},q_n)\cdot \operatorname{hor}\,\delta q_n.
\label{e-preDR}
\end{eqnarray}
Assuming that $\delta\mathbf{x}$ vanishes
at the endpoints, $\operatorname{hor}\,\delta q_0 = 0$, and
$\operatorname{hor}\,\delta q_1 = 0$ and consequently, the boundary terms vanish and we obtain the {\bfi reduced discrete variational principle},
\begin{equation} \label{e-drvar}
\delta \sum_{k=0}^{n-1} \hat{L}_d(x_k,x_{k+1})=\sum_{k=0}^{n-1}
\hat{\mathcal{A}}(x_k,x_{k+1})\cdot (\delta x_k,\delta x_{k+1}).
\end{equation}
In an analogous fashion to rewriting $D\hat L_d (x_k, x_{k+1})\cdot(\delta x_k, \delta x_{k+1})$ as $D_1 \hat L_d(x_k, x_{k+1})\cdot\delta x_k + D_2 \hat L_d(x_k, x_{k+1})\cdot\delta x_{k+1}$, we do the same by for the $\hat{\mathcal{A}}$ term by defining,
\[
\hat{\mathcal{A}}(x_0,x_1)\cdot (\delta x_0,\delta x_1)=
\hat{\mathcal{A}}_1(x_0,x_1)\cdot \delta x_0 +
    \hat{\mathcal{A}}_2(x_0,x_1)\cdot \delta x_1 .
\]
Then, equating terms involving $\delta x_k$ on the left hand side of
\eqref{e-drvar} to
the corresponding terms on the right, we get the {\bfi discrete
Routh (DR) equations} giving dynamics on $S\times S$:
\begin{equation}
\label{e-dreq}
D_2 \hat{L}_d (x_{k-1},x_k)+D_1 \hat{L}_d(x_k,x_{k+1}) =
\hat{\mathcal{A}}_2(x_{k-1},x_k) + \hat{\mathcal{A}}_1(x_k,x_{k+1}).
\end{equation}
Note that these equations depend on the value of momentum $\mu$.
Thus, if $\mathbf{q}$ is a discrete curve satisfying the
discrete Euler-Lagrange equations, the curve $\mathbf{x}$ obtained by
projecting $\mathbf{q}$ down to $S$ satisfies the DR equations \eqref{e-dreq}.

Now consider the converse, the discrete reconstruction
procedure: Given a discrete curve $\mathbf{x}$ on $S$ that satisfies the
DR equations, is $\mathbf{x}$ the projection of a discrete curve
$\mathbf{q}$ on $Q$ that satisfies the DEL equations?

Let the pair $(q_0,q_1)$ be a lift of $(x_0,x_1)$ such that
$\mathbf{J}_d (q_0,q_1)=\mu$. Let $\mathbf{q}=\{q_0,\ldots,q_n\}$ be
the solution
of the DEL equations with initial condition $(q_0,q_1)$.
Note that $\mathbf{q}$ has momentum $\mu$.
Let $\mathbf{x}' = \{x_0',\ldots,x_n'\}$ be the curve on $S$ obtained by
projecting $\mathbf{q}$. By our arguments above, $\mathbf{x}'$ solves the
DR equations. However $\mathbf{x}'$ has the
initial condition $(x_0,x_1)$, which is the same as the initial condition of
$\mathbf{x}$. By uniqueness of the solutions of the DR
equations, $\mathbf{x}'=\mathbf{x}$.
Thus $\mathbf{x}$ is the projection of a solution $\mathbf{q}$ of the DEL
equations with momentum $\mu$. Also, for a given initial
condition $q_0$, there is a unique lift of $\mathbf{x}$ to a
curve with momentum
$\mu$. Such a lift can be constructed using the method described
in Lemma~\ref{reconstruct-lemma}. Thus, lifting $\mathbf{x}$ to a curve with
momentum $\mu$ yields a solution of the discrete Euler-Lagrange equations,
which projects down to $\mathbf{x}$.

We summarize the results of this section in the following Theorem.
\begin{theorem}\label{t-drr}
Let $\mathbf{x}$ is a discrete curve on $S$, and let $\mathbf{q}$
be a discrete curve on $Q$ with momentum $\mu$ that is obtained
by lifting $\mathbf{x}$. Then the following are equivalent.
\begin{enumerate}
\item[{\rm 1.}]
$\mathbf{q}$ solves the DEL equations.
\item[{\rm 2.}]
$\mathbf{q}$ is a solution of the discrete Hamilton's variational principle
\[
\fl\qquad\delta \sum_{k=0}^{n-1}L_d(q_k,q_{k+1}) = 0
\]
for all variations $\delta\mathbf{q}$ of $\mathbf{q}$ that vanish at the
endpoints.
\item[{\rm 3.}]
$\mathbf{x}$ solves the DR equations
\[
\fl\qquad D_2 \hat{L}_d (x_{k-1},x_k)+D_1 \hat{L}_d(x_k,x_{k+1}) =
\hat{\mathcal{A}}_2(x_{k-1},x_k) + \hat{\mathcal{A}}_1(x_k,x_{k+1}).
\]
\item[{\rm 4.}]
$\mathbf{x}$ is a solution of the reduced variational principle
\[
\fl\qquad \delta \sum_{k=0}^{n-1} \hat{L}_d(x_k,x_{k+1}) =
\sum_{k=0}^{n-1} \hat{\mathcal{A}}(x_k,x_{k+1})\cdot
(\delta x_k,\delta x_{k+1})
\]
for all variations $\delta\mathbf{x}$ of $\mathbf{x}$ that vanish at the
endpoints.
\end{enumerate}
\end{theorem}

Note that for smooth group actions, the order of accuracy will be 
equal for the reduced and unreduced algorithm.

\subsection{Preservation of the Reduced Discrete Symplectic
Form}\label{s-symp} The DR equations define a discrete flow map
$\hat{F}_k:S\times S
\to S\times S$. We already know that the flow of the DEL equations preserves
the symplectic form $\Omega_{L_d}$ on $Q\times Q$. In this section
we show that the
reduced flow $\hat{F}_k$  preserves a {\em reduced} symplectic form
$\Omega_{\mu,d}$ on
$S\times S$, and that this reduced symplectic form is obtained by
restricting $\Omega_{L_d}$ to $\mathbf{J}_d ^{-1}(\mu)$ and then dropping to
$S\times S$. In other words,
$\pi_{\mu,d}^* \Omega_{\mu,d} = \iota_{\mu,d}^* \Omega_{L_d}.$
The continuous analogue of this equation is
$\pi_\mu^*\Omega_\mu = i_\mu^*\Omega_Q.$

Recall from continuous reduction theory on the Hamiltonian side that in the case of cotangent bundles, the projection $\pi_\mu:J^{-1}(\mu)\rightarrow T^*S$ can be defined as follows:
If $\alpha_q\in J^{-1}(\mu)$, then the {\bfi momentum shift}\index{momentum shift} $\alpha_q -
\mathfrak{A}_\mu (q)$ annihilates all vertical tangent vectors at $q\in
Q$, as shown by the following calculation:
\[
\langle\alpha_q - \mathfrak{A}_\mu (q), \xi_Q(q)\rangle =
J(\alpha_q)\cdot\xi - \langle\mu,\xi\rangle =
\langle\mu,\xi\rangle - \langle\mu,\xi\rangle = 0 .
\]
Thus, $\alpha_q - \mathfrak{A}_\mu (q)$ induces an element of $T_x^*
S$ and $\pi_\mu (\alpha_q)$ is defined to be this element.

Let $\mathbb{F}':\mathbf{J}_d^{-1}(\mu)\to J^{-1}(\mu)$ be the restriction
of $\mathbb{F}L_d$ to $\mathbf{J}_d^{-1}(\mu)$. Thus $\mathbb{F}'\circ \iota_\mu =\iota_{\mu,d}\circ\mathbb{F}L_d$, where
$\iota_\mu:J^{-1}(\mu)\to T^*Q$ and $\iota_{\mu,d}:\mathbf{J}_d^{-1}(\mu)\to Q\times Q$
are inclusions.
Define the map $\hat{\mathbb{F}}:S\times S\to T^*S$ by
$\hat{\mathbb{F}}(x_0,x_1)=D_2\hat{L}_d(x_0,x_1)-\hat{\mathcal{A}}_2(x_0,x_1).$ By equation~\eqref{e-crux} and Lemma~\ref{l-encapsulate-legality}, this map is well-defined.
The map $\hat{\mathbb{F}}$ will play the role of a reduced discrete
Legendre transform, and in contrast to the continuous theory, the momentum shift appears in the map $\hat{\mathbb{F}}$, as opposed to the projection $\pi_{\mu,d}$.
As in the continuous theory, $\Omega_{\mu,d}$ does involve magnetic terms.

\begin{lemma}\label{l-commute}
The following diagram commutes.
\[\xymatrix{
\mathbf{J}_d ^{-1}(\mu) \ar[r]^{\mathbb{F}'} \ar[d]_{\pi_{\mu,d}} &
J^{-1}(\mu) \ar[d]^{\pi_\mu}\\
S\times S \ar[r]^{\hat{\mathbb{F}}} & T^*S \\
}\]
\end{lemma}

\begin{proof}
Let $(q_0,q_1)\in \mathbf{J}_d^{-1}(\mu)$. Thus $D_2L_d(q_0,q_1)\in J^{-1}(\mu)$, and
\[
\pi_\mu (\mathbb{F}'(q_0,q_1))=\pi_\mu(D_2L_d(q_0,q_1)).\]

As noted above, $\pi_\mu(D_2L_d(q_0,q_1))$
is the element of $T^*_{x_1}S$ determined
by $(D_2L_d(q_0,q_1) - \mathfrak{A}_\mu(q_1))$.
For $\delta q_1\in T_{q_1}Q$, we have 
\begin{eqnarray*}
\fl\qquad\langle D_2L_d(q_0,q_1) - \mathfrak{A}_\mu(q_1),\delta q_1\rangle &= (DL_d - \mathcal{A})(q_0,q_1)\cdot (0,\delta q_1) \\
\fl\qquad &= (D\hat{L}_d - \hat{\mathcal{A}})(x_0,x_1)\cdot (0,\delta x_1)\\
\fl\qquad &=D_2\hat{L}_d(x_0,x_1)\cdot \delta x_1-\hat{\mathcal{A}}_2 (x_0,x_1)
\cdot\delta x_1\, ,
\end{eqnarray*}
where in the second equality, we used Lemma~\ref{l-encapsulate-legality}. Thus,
\[
\pi_\mu(D_2L_d(q_0,q_1))=D_2\hat{L}_d(x_0,x_1)-\hat{\mathcal{A}}_2(x_0,x_1),
\]
which means $\hat{\mathbb{F}}\circ \pi_{\mu,d} =
\pi_\mu \circ\mathbb{F}'$.
\end{proof}

\begin{theorem}\label{t-symp}
The flow of the DR equations preserves the symplectic form
\[
\Omega_{\mu,d}=\hat{\mathbb{F}}^*(\Omega_S - \pi_{T^*S,S}^*\beta_\mu),
\]
where $\beta_\mu$ is the $2$-form on $S$ obtained by dropping $\mathbf{d}\mathfrak{A}_\mu$.

Furthermore,
$\Omega_{\mu,d}$ can be obtained by dropping to $S\times S$ the restriction
of $\Omega_{L_d}$ to $\mathbf{J}_d ^{-1}(\mu)$. In other words,
$\pi_{\mu,d}^*\Omega_{\mu,d}=\iota_{\mu,d}^*\Omega_{L_d} .$
\end{theorem}

\begin{proof}
We give an outline of the steps involved; the details are routine to fill in.  The strategy is to first show that the
restriction to $\mathbf{J}_d ^{-1}(\mu)$ of the symplectic form $\Omega_{L_d}$
drops to a 2-form $\Omega_{\mu,d}$ on $S\times S$.  The fact that the
discrete flow on $Q\times Q$ preserves the symplectic form
$\Omega_{L_d}$ is then used to show that the reduced flow preserves
$\Omega_{\mu,d}$. The steps involved are as follows.
\begin{enumerate}
\item
Consider the 1-form $\Theta_{L_d}$ on $Q\times Q$ defined by
$\Theta_{L_d}(q_0,q_1)\cdot (\delta q_0,\delta q_1) = D_2L_d(q_0,q_1)\cdot \delta q_1$. The 1-form $\Theta_{L_d}$ is $G$-invariant, and thus the Lie derivative $\mathcal{L}_{\xi_{Q\times Q}}\Theta_{L_d}$ is zero.

\item
Since $\Omega_{L_d}=-\mathbf{d}\Theta_{L_d}$,
$\Omega_{L_d}$ is $G$-invariant. If $\iota_{\mu,d} : \mathbf{J}_d
^{-1}(\mu)\to Q\times Q$
is the inclusion, $\Theta_{L_d}'=\iota_{\mu,d}^*\Theta_{L_d}$ and
$\Omega_{L_d}'=\iota_{\mu,d}^*\Omega_{L_d}$ are the restrictions of
$\Theta_{L_d}$ and $\Omega_{L_d}$ respectively to $\mathbf{J}_d ^{-1}(\mu)$.
One checks that
$\Theta_{L_d}'$ and $\Omega_{L_d}'$ are invariant under the action of
$G$ on $\mathbf{J}_d ^{-1}(\mu)$.

\item
If $\xi_{\mathbf{J}_d ^{-1}(\mu)}$ is an infinitesimal generator on
$\mathbf{J}_d ^{-1}(\mu)$, then
\[\fl\quad
{\xi_{\mathbf{J}_d ^{-1}(\mu)}} \intprod\;\, \Omega_{L_d}'=
-{\xi_{\mathbf{J}_d ^{-1}(\mu)}} \intprod\;\, \mathbf{d}
\Theta_{L_d}'=-\mathcal{L}_{\xi_{\mathbf{J}_d ^{-1}(\mu)}}\Theta_{L_d}' +
\mathbf{d} {\xi_{\mathbf{J}_d ^{-1}(\mu)}} \intprod\;\, \Theta_{L_d}' =0.
\]
This follows from the $G$-invariance of $\Theta_{L_d}'$, and the
fact that $\Theta_{L_d}'\cdot \xi_{\mathbf{J}_d ^{-1}(\mu)} =\langle
\mu,\xi\rangle$.

\item
By steps 2 and 3, the form $\Omega_{L_d}'$ drops to
a reduced form $\Omega_{\mu,d}$ on $\mathbf{J}_d ^{-1}(\mu)/G\approx
S\times S$.
Thus, if $\pi_{\mu,d}: \mathbf{J}_d ^{-1}(\mu)\to S\times S$ is the
projection, then
$\pi_{\mu,d}^*\Omega_{\mu,d}=\Omega_{L_d}'$. Note that the closure of
$\Omega_{\mu,d}$ follows from the fact that $\Omega_{L_d}'$ is closed,
which in turn follows from the closure of $\Omega_{L_d}$ and the relation
$\Omega_{L_d}'=\iota_{\mu,d}^*\Omega_{L_d}$.

\item
If $F_k:Q\times Q\to Q\times Q$ is the flow of the DEL equations,
let $F_k'$ be the restriction of this flow to $\mathbf{J}_d
^{-1}(\mu)$. We know that
$F_k'$ drops to the flow $\hat{F}_k$ of the DR equations
on $S\times S$. Since $F_k$ preserves $\Omega_{L_d}$,
$F_k'$ preserves $\Omega_{L_d}'$. Using this, it can be shown that
$\hat{F}_k$ preserves $\Omega_{\mu,d}$. Note that it is sufficient
to show that $\pi_{\mu,d}^* (\hat{F}_k^* \Omega_{\mu,d}) = \pi_{\mu,d}^*
\Omega_{\mu,d}$.

\item
It now remains to compute a formula for the reduced form $\Omega_{\mu,d}$.
Using Lemma \ref{l-commute}, it follows that
\begin{eqnarray*}
\fl\qquad \pi_{\mu,d}^*\Omega_{\mu,d} &= \iota_{\mu,d}^*\Omega_{L_d}
= \iota_{\mu,d}^*\mathbb{F}L_d^*\Omega_Q = (\mathbb{F}')^*i_\mu^*\Omega_Q\\
\fl\qquad &=(\mathbb{F}')^*\pi_\mu^*(\Omega_S-\pi_{T^*S,S}^*\beta_\mu)=\pi_{\mu,d}^*\hat{\mathbb{F}}^*(\Omega_S-\pi_{T^*S,S}^*\beta_\mu).
\end{eqnarray*}
Thus $\pi_{\mu,d}^*\Omega_{\mu,d}=\pi_{\mu,d}^*\hat{\mathbb{F}}^*
(\Omega_S-\pi_{T^*S,S}^*\beta_\mu)$, from which it follows that
$\Omega_{\mu,d}=\hat{\mathbb{F}}^*(\Omega_S-\pi_{T^*S,S}^*\beta_\mu)$.
Incidentally, this expression shows that $\Omega_{\mu,d}$ is nondegenerate
provided the map $\hat{\mathbb{F}}=D_2\hat{L}_d-\hat{\mathcal{A}}_2$
is a local diffeomorphism.
\end{enumerate}
\end{proof}

One can alternatively prove symplecticity of the reduced flow directly from the reduced variational principle---see \S2.3.4 of \cite{Leok2004}.

\subsection{Relating Discrete and Continuous Reduction}
\label{s-disCts}

As shown in \cite{MaWe2001}, if the discrete
Lagrangian $L_d$ approximates the Jacobi solution of the
Hamilton-Jacobi equation, then the DEL equations give us an
integration scheme for the EL equations. In our commutative diagrams
we will denote the relationship between the EL and DEL equations by a
dashed arrow as follows:
\[
\xymatrix{
     (TQ,EL) \ar@{-->}[r] & (Q\times Q, DEL).}
\]
This arrow  can thus be read as ``the corresponding discretization''. By the continuous and discrete Noether theorems, we can restrict the flow of the EL and DEL equations to
$J_L^{-1}(\mu)$ and $\mathbf{J}_d ^{-1}(\mu)$ respectively. The flow on $J_L^{-1}(\mu)$ induces a reduced flow on
$J_L^{-1}(\mu)/G\approx TS$, which is the flow of the Routh
equations.  Similarly, the discrete flow on $\mathbf{J}_d ^{-1}(\mu)$
induces a reduced discrete flow on $\mathbf{J}_d ^{-1}(\mu)/G\approx S\times
S$, which is the flow of the discrete Routh equations.  Since the
DEL equations give us an integration algorithm for the EL
equations, it follows that the DR equations give us an integration
algorithm for the Routh equations.

To numerically integrate the Routh equations, we have two options:
\begin{enumerate}
\item First solve the DEL equations to yield a discrete trajectory on
     $Q$, which can then be projected to a discrete trajectory on $S$.
\item Solve the DR equations to directly obtain a discrete trajectory
     on $Q$.
\end{enumerate}
Either approach yields the same result, which we express by the commutative diagram:
\begin{equation}\label{e-cd-crDr} \xymatrix{
       (J_L^{-1}(\mu),EL) \ar@{-->}[r] \ar[d]_{\pi_{\mu,L}}
       & (\mathbf{J}_d ^{-1}(\mu),DEL) \ar[d]^{\pi_{\mu,d}}\\
       (TS, R) \ar@{-->}[r] & (S\times S,DR)
}
\end{equation}

\section{Reduction of the symplectic Runge-Kutta algorithm.}\label{s-srkRed}
A well-studied class of numerical schemes for Hamiltonian and
Lagrangian systems is the symplectic partitioned Runge-Kutta (SPRK) algorithms
(see \cite{HaNoWa1993, HaWa1996} for history and details). We will adopt a local
trivialization to express the SPRK method in which the group action is addition. Given a connection $\mathfrak{A}$ on $Q$, it can be represented in local coordinates as $\mathfrak{A}(\theta,x)(\dot{\theta},\dot{x})=A(x)\dot{x}+\dot{\theta} .$ By rewriting the symplectic partitioned Runge-Kutta algorithm in terms of this local trivialization, and using the local representation of the connection, we obtain the following algorithm on $T^*S$:
\numparts
\begin{eqnarray}\label{e-rsrk-start}
\fl\qquad x_1&=x_0+h\sum b_j \dot{X}_j \\
\fl\qquad s_1&=s_0+h\sum_j\tilde{b}_j\dot{S}_j
+\left[ h
\sum_j \left( \tilde{b}_j \mu \frac{\partial A}{\partial x}
(X_j)\dot{X}_j\right)-\left(\mu A(x_1)-\mu A(x_0)\right)\right]\\
\fl\qquad X_i&=x_0+h\sum a_{ij} \dot{X}_j\\
\fl\qquad S_i&=s_0+h\sum_j\tilde{a}_{ij}\dot{S}_j
+\left[ h
\sum_j \left( \tilde{a}_{ij} \mu \frac{\partial A}{\partial x}
(X_j)\dot{X}_j\right)-\left(\mu A(X_i)-\mu A(x_0)\right)\right]\\
\fl\qquad S_j&=\frac{\partial \hat{R}^\mu}{\partial \dot{x}}(X_j,\dot{X}_j)\\
\fl\qquad \dot{S}_j&= \frac{\partial \hat{R}^\mu}{\partial x}(X_j,\dot{X}_j)
-i_{\dot{X}_j}\beta_\mu(X_j),\label{e-rsrk-end}
\end{eqnarray}
\endnumparts
This system of equations is called the {\bfi reduced symplectic
partitioned Runge-Kutta (RSPRK)} algorithm. A detailed derivation can be found in  \S2.5 of \cite{Leok2004}. As we obtained this system
by dropping the symplectic
partitioned Runge-Kutta algorithm from $\mathbf{J}^{-1}(\mu)$ to $T^*S$, it follows that
this algorithm preserves the reduced symplectic form
$\Omega_\mu = \Omega_S -\pi_{T^*S,S}^*\beta_\mu$ on $T^*S$.

Since the SPRK algorithm is an integration algorithm for the
Hamiltonian vector field $X_H$ on $T^{\ast}Q$, the RSPRK algorithm is
an integration
algorithm for the reduced Hamiltonian vector field $X_{H_\mu}$ on
$T^*S$. The relationship between cotangent bundle reduction and the
reduction of the SPRK algorithm can be represented
by the following commutative diagram:
\[\xymatrix{
(J^{-1}(\mu),X_H) \ar[d]_{\pi_{\mu}}\ar@{-->}[r] & (J^{-1}(\mu),SPRK)
\ar[d]^{\pi_{\mu}}\\
(T^*S, X_{H_\mu}) \ar@{-->}[r]& (T^*S,RSPRK)
}
\]
The dashed arrows here denote the corresponding discretization, as
in \eqref{e-cd-crDr}.
The SPRK algorithm can be obtained by pushing forward the DEL equations by
the discrete Legendre transform. See, for example, \cite{MaWe2001}. By Lemma \ref{l-commute}, this
implies that the RSPRK algorithm can be obtained by pushing
forward the DR equations by the reduced discrete Legendre
transform $\hat{\mathbb{F}}=D_2 \hat{L}_d - \hat{\mathcal{A}}_2$.
These relationships are shown in the following commutative
diagram:
\[\xymatrix{
(\mathbf{J}_d ^{-1}(\mu), DEL) \ar[r]^{\mathbb{F}L_d}
\ar[d]_{\pi_{\mu,d}}& (J^{-1}(\mu),SPRK) \ar[d]^{\pi_\mu} \\
(S\times S,DR) \ar[r]^{\hat{\mathbb{F}}} & (T^*S, RSPRK)
}
\]

\section{Putting Everything Together}\label{s-summary}
The relationship between Routh reduction and cotangent bundle reduction can be represented by the following commutative diagram:
\[\xymatrix{
(J^{-1}_L(\mu),EL) \ar[r]^{\mathbb{F}L} \ar[d]_{\pi_{\mu,L}} &
(J^{-1}(\mu),X_H) \ar[d]^{\pi_{\mu}}\\
(TS,R) \ar[r]^{\mathbb{F}\hat{R}^{\mu}} & (T^{*}S,X_{H_{\mu}})
}
\]
We saw in \S\ref{s-disCts} that if $L_d$ approximates the Jacobi
solution of the Hamilton-Jacobi equation, discrete and continuous Routh reduction is related by the
following diagram:
\[\xymatrix{
(J^{-1}_L(\mu),EL) \ar@{-->}[r] \ar[d]_{\pi_{\mu,L}} &
(J^{-1}_d(\mu),DEL) \ar[d]^{\pi_{\mu,d}}\\
(TS,R) \ar@{-->}[r] & (S\times S,DR)
}
\]
The dashed arrows mean that the DEL equations are an integration
algorithm for the EL equations, and that the DR equations are an
integration algorithm for the Routh equations.

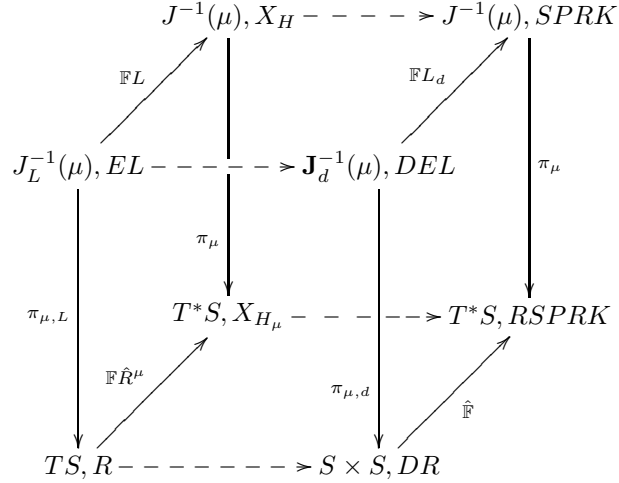
\begin{figure}
\begin{equation*}
    \xymatrix@!0@R=2cm@C=2cm{
      & J^{-1}(\mu),X_H \ar@{-->}[rr]
      \ar'[d][dd]_{\pi_{\mu}} & & J^{-1}(\mu),SPRK \ar[dd]^{\pi_{\mu}} \\
      J_L^{-1}(\mu),EL \ar@{-->}[rr] \ar[ur]^{\mathbb{F}L}
      \ar[dd]_{\pi_{\mu,L}} & & \mathbf{J}_d ^{-1}(\mu),DEL
      \ar[ur]^{\mathbb{F}L_d} \ar[dd]_(0.75){\pi_{\mu,d}} \\
      & T^*S,X_{H_{\mu}} \ar@{-->}'[r][rr] & & T^*S,RSPRK \\
      TS,R \ar@{-->}[rr] \ar[ur]^{\mathbb{F}\hat{R}^{\mu}}
      & & S \times S,DR \ar[ur]_{\hat{\mathbb{F}}}
      }
\end{equation*}
\caption{
     Complete commutative cube. Dashed arrows represent discretization
     from the continuous systems on the left face to the discrete systems
     on the right face. Vertical arrows represent reduction from the full
     systems on the top face to the reduced systems on the bottom face.
     Front and back faces represent Lagrangian and Hamiltonian
     viewpoints, respectively.
\label{fig:diagram}}
\end{figure}

Pushing forward the DEL equation using the discrete Legendre transform $\mathbb{F}L_d$ yields the SPRK algorithm on $T^{\ast}Q$,
which is an
integration algorithm for $X_H$. This is depicted by:
\[\xymatrix{
(J^{-1}_L(\mu),EL) \ar@{-->}[r] \ar[d]_{\mathbb{F}L} &
(J^{-1}_d(\mu),DEL) \ar[d]^{\mathbb{F}L_d}\\
(J^{-1}(\mu),X_H) \ar@{-->}[r] & (J^{-1}(\mu),SPRK)
}
\]
The SPRK algorithm on $J^{-1}(\mu)\subset T^{\ast}Q$ induces the RSPRK
algorithm on $J^{-1}(\mu)/G \approx T^*S$. As we saw in
\S\ref{s-srkRed}, this reduction process is related to cotangent
bundle reduction and to discrete Routh reduction as shown in the
following diagram:
\[\xymatrix{
(J^{-1}(\mu),X_H) \ar@{-->}[r] \ar[d]_{\pi_{\mu}} &
(J^{-1}(\mu),SPRK) \ar[d]^{\pi_{\mu}} &
(J^{-1}_d(\mu),DEL) \ar[l]_{\mathbb{F}L_d} \ar[d]^{\pi_{\mu,d}}\\
(T^{*}S,X_{H_{\mu}}) \ar@{-->}[r] & (T^{*}S,RSPRK) & (S\times
S,DR)\ar[l]_{\hat{\mathbb{F}}}
}
\]
Putting all the above commutative diagrams together into one
diagram, we obtain Figure \ref{fig:diagram}.

\section{Example: $J_2$ Satellite Dynamics}
\label{s-examples}

\subsection{Configuration Space and Lagrangian}

An illustrative and important example of a system with an abelian
symmetry group is that of a single satellite in orbit about an oblate
Earth. The general aspects and background for this problem are discussed
in \cite{PrCo1993}, and some interesting aspects of the geometry
underlying it, are discussed in \cite{ChMa2003}.

The configuration manifold $Q$ is $\mathbb{R}^3$, and the Lagrangian
is
\[
L(q,\dot{q}) = \frac{1}{2} M_s \|\dot{q}\|^2 - M_s V(q),
\]
where $M_s$ is the mass of the satellite and $V : \mathbb{R}^3 \to
\mathbb{R}$ is the gravitational potential due to the earth truncated
at the first term in the expansion in the ellipticity
\[
V(q) = \frac{G M_e}{\|q\|} + \frac{G M_e R_e^2 J_2}{\|q\|^3}
\left(\frac{3}{2} \frac{(q^3)^2}{\|q\|^2} - \frac{1}{2}\right).
\]
Here, $G$ is the gravitational constant, $M_e$ is the mass of the
Earth, $R_e$ is the radius of the earth, $J_2$ is a small
non-dimensional parameter describing the degree of ellipticity, and
$q^3$ is the third component of $q$. In non-dimensional coordinates,
\begin{equation}
     \label{eqn:sat_L}
     L(q,\dot{q}) = \frac{1}{2}\|\dot{q}\|^2
     - \left[ \frac{1}{\|q\|} + \frac{J_2}{\|q\|^3}
       \left(\frac{3}{2} \frac{(q^3)^2}{\|q\|^2} - \frac{1}{2}\right)
     \right].
\end{equation}
This corresponds to choosing space and time coordinates in which the
radius of the Earth is $1$ and the period of orbit at zero altitude is
$2\pi$ when $J_2 = 0$ (spherical earth).

\subsection{Symmetry Action}

The symmetry of interest to us is that of rotation about the vertical
($q^3$) axis, so the symmetry group is the unit circle $S^1$. Using
cylindrical coordinates $q = (r,\theta,z)$ for the configuration, the
symmetry action is $\phi : (r,\theta,z) \mapsto (r,\theta + \phi,z)$.
Since $\|q\|$, $\|\dot{q}\|$, and $q^3 = z$ are all invariant under
this transformation, so too is the Lagrangian.

This action is clearly not free on all of $Q = \mathbb{R}^3$, as the
$z$-axis is invariant for all group elements. This is not a serious
obstacle as the lifted action is free on $T(Q \setminus
(0,0,0))$ and this is enough to permit the application of the intrinsic
Routh reduction theory.
Alternatively,
one can simply take $Q = \mathbb{R}^3 \setminus \{(0,0,z)\mid z \in
\mathbb{R}\}$ and then the theory literally applies.

The shape space $S = Q / G$ is thus the half-plane $S = \mathbb{R}^+
\times \mathbb{R}$ and we will take coordinates $(r,z)$ on $S$. In
doing so, we are implicitly defining a global diffeomorphism $S \times
G \to Q$ given by $((r,z),\theta) \mapsto (r,\theta,z)$.

The Lie algebra $\mathfrak{g}$ for $G = S^1$ is the real line
$\mathfrak{g} = \mathbb{R}$, and we will identify the dual with the
real line itself $\mathfrak{g}^* \cong \mathbb{R}$. For a Lie algebra
element $\xi \in \mathfrak{g}$, the corresponding infinitesimal
generator is given by
$\xi_Q : (r,\theta,z) \mapsto ((r,\theta,z),(0,\xi,0)).$
The Lagrange momentum map $\mathbf{J}_L : TQ \to \mathfrak{g}^*$ is
given by
$\mathbf{J}_L(v_q)\cdot \xi =\langle\mathbb{F}L(v_q),\xi_Q(q)\rangle,$
which in our case is a scalar quantity, the vertical component of the
standard angular momentum:
$J_L((r,\theta,z),(\dot{r},\dot{\theta},\dot{z})) =
r^2\dot{\theta}.$

Consider the Euclidean metric on $\mathbb{R}^3$, which corresponds to
the kinetic energy norm in the Lagrangian. From this metric we define
the mechanical connection $\mathfrak{A} : TQ \to \mathfrak{g}$ given
by $\mathfrak{A}((r,\theta,z),(\dot{r},\dot{\theta},\dot{z})) =
\dot{\theta}$. The $1$-form $\mathfrak{A}_{\mu}$ on $Q$ is thus given
by $\mathfrak{A}_{\mu} = \mu \mathbf{d}\theta$. The exterior
derivative of this expression gives $\mathbf{d}\mathfrak{A}_{\mu} =
\mu \mathbf{d}^2\theta = 0$, and so the reduced $2$-form is
$\beta_{\mu} = 0$.

\subsection{Equations of Motion}

The Euler-Lagrange equations for the Lagrangian
\eqref{eqn:sat_L} gives the equations of motion,
\[
\ddot{q} = - \nabla_q \left[ \frac{1}{\|q\|} + \frac{J_2}{\|q\|^3}
     \left(\frac{3}{2} \frac{(q^3)^2}{\|q\|^2} - \frac{1}{2}\right)
\right].
\]
To calculate the reduced equations, we begin by calculating the
Routhian
\begin{eqnarray*}
     R^{\mu}(r,\theta,z,\dot{r},\dot{\theta},\dot{z}) &=
     L(r,\theta,z,\dot{r},\dot{\theta},\dot{z})
     - \mathfrak{A}_{\mu}(r,\theta,z) \cdot (\dot{r},\dot{\theta},\dot{z}) \\
     &= \frac{1}{2}\|(\dot{r},\dot{\theta},\dot{z})\|^2 - \left[
       \frac{1}{r} + \frac{J_2}{r^3} \left(\frac{3}{2} \frac{z^2}{r^2} -
         \frac{1}{2}\right) \right] - \mu\dot{\theta}.
\end{eqnarray*}
We choose a fixed value $\mu$ of the momentum and restrict
ourselves to the space $J_L^{-1}(\mu)$, on which $\dot{\theta} = \mu$.
The reduced Routhian $\hat{R}^{\mu} : TS \to \mathbb{R}$ is
the restricted Routhian dropped to the tangent bundle of the shape
space. In coordinates this is
\begin{equation*}
     \hat{R}^{\mu}(r,z,\dot{r},\dot{z})
     = \frac{1}{2}\|(\dot{r},\dot{z})\|^2
     - \left[ \frac{1}{r} + \frac{J_2}{r^3}
       \left(\frac{3}{2} \frac{z^2}{r^2} - \frac{1}{2}\right)
     \right] - \frac{1}{2}\mu^2.
\end{equation*}
Recalling that $\beta_{\mu} = 0$, the Routh equations
can be evaluated to give
\[
(\ddot{r},\ddot{z}) = - \nabla_{(r,z)}
\left[ \frac{1}{r} + \frac{J_2}{r^3}
\left(\frac{3}{2} \frac{z^2}{r^2} - \frac{1}{2}\right)
\right],
\]
which describes the motion on the shape space.

To recover the unreduced Euler-Lagrange equations from the Routh
equations one uses the procedure of reconstruction. This is covered in
detail in \cite{MaMoRa1990, Marsden1992, MaRaSc2000}.

\subsection{Discrete Lagrangian System}

We discretize this system with a high-order discrete Lagrangian.
Recall that the pushforward discrete Lagrange map associated with this discrete Lagrangian is a symplectic partitioned Runge-Kutta method.

Given a point $(q_0,q_1) \in Q \times Q$ we will take $(q_0,p_0)$ and
$(q_1,p_1)$ to be the associated discrete Legendre transforms. As the
discrete momentum map is the pullback of the canonical momentum map,
we have that
$J_{L_d}(q_0,q_1) = (p_{\theta})_0 = (p_{\theta})_1.$
Take a fixed momentum map value $\mu$ and restrict $L_d$ to the set
$J_{L_d}^{-1}(\mu)$. Dropping this to $S \times S$ now gives the
reduced discrete Lagrangian $\hat{L}_d : S \times S \to \mathbb{R}$. More explicitly, $L_d$ depends on $(r_k, \theta_k, z_k)$ and $(r_{k+1},\theta_{k+1},z_{k+1})$, but group invariance implies that the group variables only enters in the combination $(\theta_{k+1}-\theta_k)$. The condition $\mathbf{J}_d(q_k,q_{k+1})=\mu$ can be inverted to eliminate the dependence of $L_d$ on $(\theta_{k+1}-\theta_k)$, and $L_d$ can be expressed in terms of $\mu$, $(r_k, z_k)$, and $(r_{k+1}, z_{k+1})$, which yields $\hat{L}_d$.

As discussed in \S\ref{s-srkRed}, the fact that we have taken
coordinates in which the group action is addition in $\theta$ means
that the pushforward discrete Lagrange map associated with the reduced
discrete Lagrangian is the reduced method given by (\ref{e-rsrk-start}-\ref{e-rsrk-end}). In
fact, as the mechanical connection has $A(r,z) = 0$ and $\beta_{\mu} =
0$, the pushforward discrete Lagrange map is exactly a partitioned
Runge-Kutta method with Hamiltonian equal to the reduced Routhian. These are generically related by a momentum shift, rather than being equal.

Given a trajectory of the reduced discrete system, we can reconstruct
the unreduced discrete trajectory by solving for $\theta$. Correspondingly, a trajectory
of the unreduced discrete system can be projected onto shape space
to give a trajectory of the reduced discrete system.

\subsection{Example Trajectories}
We compute the unreduced trajectories using the fourth-order SPRK algorithm, and the reduced trajectories using the corresponding fourth-order RSPRK algorithm.

\medskip 

\begin{figure}[h]
     \begin{center}
       \includegraphics[width=10cm]{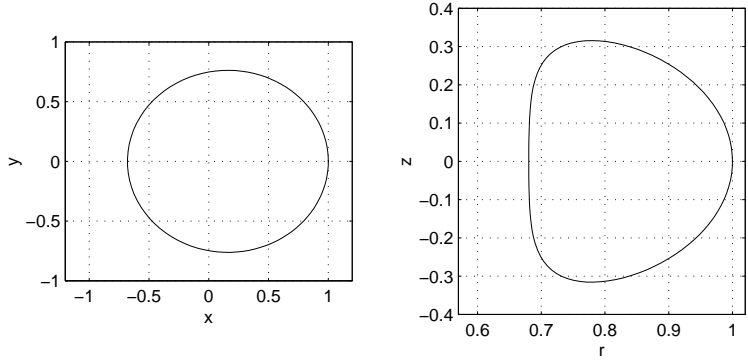}
     \end{center}
     \caption{Unreduced (left) and reduced (right) views of an inclined
       elliptic trajectory for the continuous time system with $J_2 = 0$
       (spherical Earth).
       \label{fig:spherical}}
\end{figure}

\noindent {\bf Solutions of the Spherical Earth System.}
Consider initially the system with $J_2 = 0$. This corresponds to the
case of a spherical Earth, and so the equations reduce to the standard
Kepler problem. A slightly inclined circular trajectory is shown in Figure
\ref{fig:spherical}, in both the unreduced and reduced pictures. Note
that the graph of the reduced trajectory is a quadratic, as $\|q\| =
\sqrt{r^2+z^2}$ is a constant.

We will now investigate the effect of two different perturbations to
the system, one due to taking non-zero $J_2$ and the other due to the
numerical discretization.

\medskip 

\begin{figure}[h]
     \begin{center}
       \includegraphics[width=10cm]{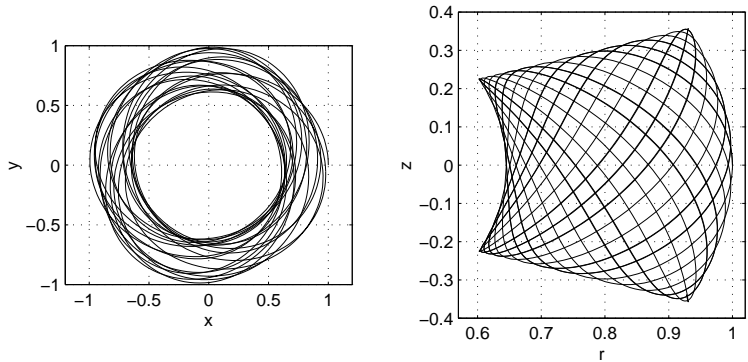}
     \end{center}
     \caption{Unreduced (left) and reduced (right) views of an inclined
       elliptic trajectory for the continuous time system with $J_2 =
       0.05$. Observe that the non-spherical terms introduce precession
       of the near-elliptic orbit in the symmetry direction.
       \label{fig:j2}}
\end{figure}

\noindent {\bf The $J_2$ effect.}
Taking $J_2 = 0.05$ (which is close the actual value for the Earth),
the system becomes near-integrable and experiences breakup of the KAM
tori. This can be seen in Figure \ref{fig:j2}, where the same initial
condition is used as in Figure \ref{fig:spherical}.

Due to the fact that the reduced trajectory is no longer a simple
curve, there is a geometric-phase-like effect which causes precession
of the orbit. This precession can be seen in the thickening of the
unreduced trajectory.

\medskip 

\begin{figure}[h]
     \begin{center}
       \includegraphics[width=10cm]{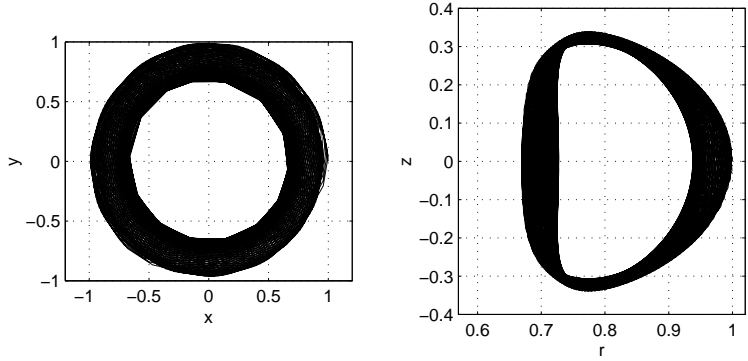}
     \end{center}
     \caption{Unreduced (left) and reduced (right) views of an inclined
       trajectory of the discrete system with step-size $h = 0.3$ and
       $J_2 = 0$. The initial condition is the same as that used in
       figure \ref{fig:spherical}. The numerically introduced precession
       means that the unreduced picture looks similar to that of figure
       \ref{fig:j2} with non-zero $J_2$, whereas only by considering the
       reduced picture we can see the correct resemblance to the $J_2 =
       0$ case of figure \ref{fig:spherical}.
       \label{fig:spherical_numerical}}
\end{figure}

\noindent {\bf Solutions of the Discrete System for a Spherical Earth.}
We now consider the discrete system with $J_2 = 0$, for the second
order Gauss-Legendre discrete Lagrangian with stepsize of $h = 0.3$.
The trajectory with the same initial condition as above is given in
Figure \ref{fig:spherical_numerical}.

As can be seen from the reduced trajectory, the discretization has
caused a similar breakup of the periodic orbit as was produced by the
non-zero $J_2$. This induces precession
of the orbit in the unreduced trajectory, in a way which is difficult
to distinguish from the perturbation above due to non-zero $J_2$ if
only the unreduced picture is considered. If the reduced pictures are
consulted, however, then it is immediately clear that the system is
much closer to the continuous time system with $J_2 = 0$ than to the
system with non-zero $J_2$.

\medskip 

\begin{figure}
     \begin{center}
       \includegraphics[width=10cm]{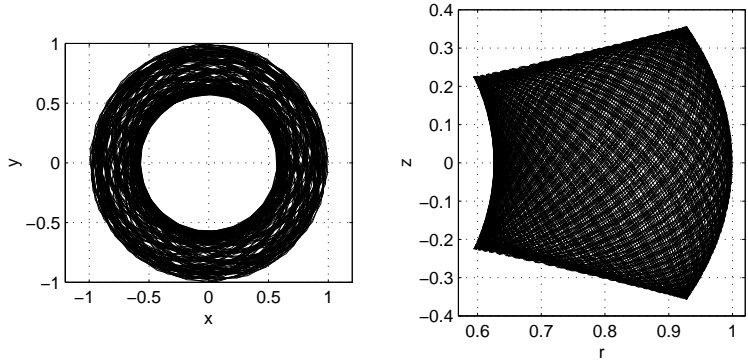}
     \end{center}
     \caption{Unreduced (left) and reduced (right) views of an inclined
       trajectory of the discrete system with step-size $h = 0.3$ and
       $J_2 = 0.05$. The initial condition is the same as that used in
       figure \ref{fig:j2}. The unreduced picture is similar to that of
       both figures \ref{fig:j2} and \ref{fig:spherical_numerical}. By
       considering the reduced picture the correct resemblance to
       \ref{fig:j2}.
       \label{fig:j2_numerical}}
\end{figure}

\noindent {\bf Solutions of the Discrete System with $J_2$ Effect.}
Finally, we consider the discrete system with non-zero $J_2 = 0.05$.
The resulting trajectory is shown in Figure \ref{fig:j2_numerical},
and it is clearly not easy to determine from the unreduced picture
whether the precession is due to the $J_2$ perturbation, the
discretization, or some combination of the two.

Taking the reduced trajectories, however, immediately shows that this
discrete time system is structurally much closer to the non-zero $J_2$
system than to the original $J_2 = 0$ system. This confusion arises
because both the $J_2$ term and the discretization introduce
perturbations which act in the symmetry direction.

While this system is sufficiently simple that one can run simulations
with such small timesteps that the discretization artefacts become
negligible, this is certainly not generally possible. This example
demonstrates how knowledge of the geometry of the system can be
important in understanding the discretization process, and how this
can give insight into the behavior of numerical simulations. In
particular, understanding how the discretization interacts with the
symmetry action is extremely important.

\subsection{Coordinate Systems}

In this example we have chosen cylindrical coordinates, thus making
the group action addition in $\theta$. One can always do this, as an
abelian Lie group is isomorphic to a product of copies of $\mathbb{R}$
and $S^1$, but it may sometimes be preferable to work in coordinates
in which the group action is not addition. For example, cartesian
coordinates in the present example. Reasons for choosing a different coordinate system might include ease of computation, or simplicity of the expressions.

If we adopt a coordinate system wherein the group action is not expressed in terms of addition, the RSPRK method is not applicable, but we can still apply the Discrete Routh equations to obtain an integration scheme on $S\times S$. The push forward of this under $\hat{\mathbb{F}}$
yields an integration scheme on $T^*S$. The trajectories on the
shape space that we obtain in this manner could be different from
those we would get with the RSPRK method. However, in both cases we
would have conservation of symplectic structure, momentum, and the
order of accuracy would be the same. One could choose whichever
approach is cheaper and easier.

\section{Example: Double Spherical Pendulum}
\label{s-examples-dsp}

\subsection{Configuration Space and Lagrangian}

We consider the example of the double spherical pendulum which has
a non-trivial magnetic term and constraints. The configuration
manifold $Q$ is $S^2\times S^2$, and the embedding linear space
$V$ is $\mathbb{R}^3\times \mathbb{R}^3$. The position vectors of each pendulum with respect to the pivot
point are denoted by $\mathbf{q}_1$ and $\mathbf{q}_2$. These
vectors are constrained to have lengths $l_1$ and $l_2$
respectively, and the pendula masses are denoted by $m_1$ and
$m_2$. The Lagrangian is
\[
\fl\qquad L(\mathbf{q}_1, \mathbf{q}_2, \dot{\mathbf{q}_1},
\dot{\mathbf{q}_2}) = \frac{1}{2} m_1 \|\dot{\mathbf{q}_1}\|^2 +
\frac{1}{2} m_2 \|\dot{\mathbf{q}_1}+\dot{\mathbf{q}_2}\|^2- m_1 g
\mathbf{q}_1 \cdot \mathbf{k} -m_2 g (\mathbf{q}_1 +
\mathbf{q}_2)\cdot \mathbf{k},
\]
where $g$ is the gravitational constant, and $\mathbf{k}$ is the
unit vector in the $z$ direction. The constraint function
$c:V\rightarrow \mathbb{R}^2$ is given by,
$c(\mathbf{q}_1,\mathbf{q}_2)=(\|\mathbf{q}_1\|-l_1,\|\mathbf{q}_2\|-l_2).$
Using cylindrical coordinates $\mathbf{q}_i = (r_i,\theta_i,z_i)$,
$L$ becomes
\begin{eqnarray*}
\fl\qquad L(q,\dot{q})=\frac{1}{2}m_{1}\left(
\dot{r}_{1}^{2}+r_{1}^{2}\dot{\theta}_{1}^{2}+%
\dot{z}_{1}^{2}\right) +\frac{1}{2}m_{2}\left\{ \dot{r}_{1}^{2}+r_{1}^{2}%
\dot{\theta}_{1}^{2}+\dot{r}_{2}^{2}+r_{2}^{2}\dot{\theta}_{2}^{2}\right.\\
\fl\qquad\qquad\qquad\left. +2\left(
\dot{r}_{1}\dot{r}_{2}+r_{1}r_{2}\dot{\theta}_{1}\dot{\theta%
}_{2}\right) \cos \varphi +2\left( r_{1}\dot{r}_{2}\dot{\theta}_{1}-r_{2}%
\dot{r}_{1}\dot{\theta}_{2}\right) \sin \varphi +\left( \dot{z}_{1}+\dot{z}%
_{2}\right) ^{2}\right\}\\
\fl\qquad\qquad\qquad-m_{1}gz_{1}-m_{2}g\left( z_{1}+z_{2}\right),
\end{eqnarray*}
where $\varphi=\theta_2-\theta_1$. Furthermore, we can
automatically satisfy the constraints by performing the
substitutions,
$z_{i}=(l_{i}^{2}-r_{i}^{2})^{1/2}$, and
$\dot{z}_{i}=-r_{i}\dot{r}_{i}(\sqrt{l_{i}^{2}-r_{i}^{2}})^{-1/2}$.
\subsection{Symmetry Action}

The symmetry of interest to us is the simultaneous rotation of the
two pendula about vertical ($z$) axis, so the symmetry group is
the unit circle $S^1$. Using cylindrical coordinates $\mathbf{q}_i
= (r_i,\theta_i,z_i)$ for the configuration, the symmetry action
is $\phi : (r_i,\theta_i,z_i) \mapsto (r_i,\theta_i + \phi,z_i)$.
Since $\|\mathbf{q}_i\|$, $\|\dot{\mathbf{q}_i}\|$,
$\|\dot{\mathbf{q}_1}+\dot{\mathbf{q}_2}\|$, and
$\mathbf{q}_i\cdot \mathbf{k}$ are all invariant under this
transformation, so too is the Lagrangian.

This action is clearly not free on all of $V = \mathbb{R}^3\times
\mathbb{R}^3$, as the $z$-axis is invariant for all group
elements. However, this does not pose a problem computationally,
as long as the trajectories do not pass through the downward
hanging configuration, corresponding to $r_1=r_2=0$. To treat the
downward handing configuration properly, we would need to develop
a discrete Lagrangian analogue of the continuous theory of
singular reduction described in \cite{OrRa2004}.

We will now show the geometric structures involved in implementing the RSPRK algorithm. See \S2.10.2 of \cite{Leok2004} for a detailed derivation. The Lie algebra $\mathfrak{g}$ for $G = S^1$ is the real line
$\mathfrak{g} = \mathbb{R}$, and we will identify the dual with
the real line itself $\mathfrak{g}^* \cong \mathbb{R}$. For a Lie
algebra element $\xi \in \mathfrak{g}$, the corresponding
infinitesimal generator is given by
$ \xi_Q : (r_1,\theta_1,z_1,r_2,\theta_2,z_2) \mapsto
((r_1,\theta_1,z_1,r_2,\theta_2,z_2),(0,\xi,0,0,\xi,0)).$
As such, the momentum map is given by
\begin{eqnarray*}
\fl\qquad J_L((r_1,\theta_1,z_1,r_2,\theta_2,z_2),(\dot{r}_1,\dot{\theta}_1,\dot{z}_1,\dot{r}_2,\dot{\theta}_2,\dot{z}_2))\\
\fl\quad\qquad = \left( m_{1}+m_{2}\right)
r_{1}^{2}\dot{\theta}_{1}+m_{2}r_{2}^{2}\dot{%
\theta}_{2}+m_{2}r_{1}r_{2}\left(
\dot{\theta}_{1}+\dot{\theta}_{2}\right) \cos \varphi +\left(
r_{1}\dot{r}_{2}-r_{2}\dot{r}_{1}\right) \sin \varphi,
\end{eqnarray*}
which is simply the vertical component of the standard
angular momentum.

The locked inertia tensor is given by \cite{Marsden1992},
\begin{equation*}
\fl\qquad\mathbb{I}(\mathbf{q}_1\mathbf{q}_2) =
m_1\|\mathbf{q}_1^\bot\|^2+m_2\|(\mathbf{q}_1+\mathbf{q}_2)^\bot\|^2 = m_1 r_1^2+m_2(r_1^2+r_2^2+2r_1 r_2 \cos\varphi).
\end{equation*}
The mechanical connection as a $1$-form is given by
\begin{eqnarray*}
\fl\qquad\alpha(\mathbf{q}_1,\mathbf{q}_2) = 
[m_1
r_1^2+m_2(r_1^2+r_2^2+2r_1 r_2 \cos\varphi)]^{-1}
\left[\left( m_1+m_2\right) r_1^2\mathbf{d}\theta_1+m_2 r_2^2\mathbf{d} \theta_2\right.\\
\fl\qquad\qquad\qquad\left.+m_2 r_1 r_2
\left( \mathbf{d}\theta_1+\mathbf{d}\theta_2\right) \cos \varphi
+\left( r_1
\mathbf{d}r_2-r_2\mathbf{d}r_1\right) \sin \varphi \right].
\end{eqnarray*}
The $\mu$-component of the mechanical connection is given by,
\begin{eqnarray*}
\fl\qquad\alpha_\mu(\mathbf{q}_1,\mathbf{q}_2)  = \mu [ m_1
r_1^2+m_2(r_1^2+r_2^2+2r_1 r_2 \cos\varphi)]^{-1}\\
\fl\qquad\qquad\qquad\times\left\{\left[\left( m_1+m_2\right) r_1^2+m_2 r_1
r_2 \cos \varphi \right] \mathbf{d}\theta_1 + \left[ m_2 r_2^2+m_2
r_1 r_2 \cos \varphi \right] \mathbf{d}\theta_2\right\}.
\end{eqnarray*}
Taking the exterior derivative of this $1$-form yields a
non-trivial magnetic term on the reduced space, which drops to the quotient space to yield,
\begin{eqnarray*}
\fl\qquad\beta_\mu = \mu m_2\left[ 2\left( m_1+m_2\right) r_1 r_2+\left( m_1
r_1^2+m_2 (r_1^2+r_2^2)\right) \cos \varphi\right] \\
\fl\qquad\qquad\qquad\times\left[ m_1
r_1^2+m_2\left( r_1^2+r_2^2 + 2 r_1 r_2 \cos\varphi \right)
\right]^{-2}
\mathbf{d}\varphi\wedge(r_2\mathbf{d}r_1-r_1\mathbf{d}r_2).
\end{eqnarray*}
The local representation of the connection is given by
\[\fl\qquad A(r_1,r_2,\varphi)= m_2(m_1 r_1^2+m_2(r_1^2+r_2^2+2r_1 r_2
\cos\varphi))^{-1}
\left[\begin{array}{c}
    -r_2\sin\varphi \\ r_1\sin\varphi \\ r_2^2+r_1 r_2\cos\varphi
\end{array}\right]^\mathrm{T},
\]
and the amended potential $V_\mu$ has the form
\begin{eqnarray*}
\fl\qquad V_\mu(q) = -m_1 g(l_1^2-r_1^2)^{1/2}-m_2
g\left[(l_1^2-r_1^2)^{1/2}+(l_2^2-r_2^2)^{1/2}\right]\\
\fl\qquad\qquad\qquad+2^{-1}\mu^2  [ m_1
r_1^2+m_2(r_1^2+r_2^2+2r_1 r_2 \cos\varphi) ]^{-1}.
\end{eqnarray*}
The Routhian on the momentum level set is given by,
$R^\mu=\frac{1}{2}\|\operatorname{hor}(q,v)\|^2-V_\mu.$
Recall that $\operatorname{hor}(v_q) = v_q - \xi_Q(v_q)$, where
$\xi=\alpha(v_q)$, and $\xi_Q(v_q)=(0,\xi,0,0,\xi,0)$. Then we
obtain,
\begin{equation*}
\fl\qquad\operatorname{hor}(v_q) = v_q -
(0,\alpha(v_q),0,0,\alpha(v_q),0)= (\dot{r}_1,\dot{\theta}_1 - \alpha(v_q)
,\dot{z}_1,\dot{r}_2,\dot{\theta}_2 - \alpha(v_q) ,\dot{z}_2).
\end{equation*}
The kinetic energy metric has the form
\begin{equation*}\fl\qquad\left[ \begin{array}{cccccc}
m_1+m_2             & 0                      & 0       &
m_2\cos\varphi     & -m_2 r_2\sin\varphi    & 0\\
0                   & (m_1 + m_2) r_1^2      & 0       & m_2
r_1\sin\varphi & m_2 r_1 r_2\cos\varphi & 0\\
0                   & 0                      & m_1+m_2 & 0
& 0                      & 0\\
m_2\cos\varphi      & m_2 r_1\sin\varphi     & 0       & m_2
& 0                      & 0\\
-m_2 r_2\sin\varphi & m_2 r_1 r_2\cos\varphi & 0       & 0
& m_2 r_2^2              & 0\\
0 & 0 & 0 & 0 & 0 & m_2\\
\end{array}\right]
\end{equation*}
This together with the expression for $\operatorname{hor}(v_q)$ allows us to compute $\frac{1}{2}\|\operatorname{hor}(q,v)\|^2$, and when combined with the formula for the amended potential $V_\mu$ gives the Routhian $R^\mu$. Notice that all our expressions are in terms of the reduced variables on $TS$, so they trivially drop to yield $\hat R^\mu$. These expressions can then be directly substituted into the RSPRK algorithm in equations (\ref{e-rsrk-start}-\ref{e-rsrk-end}) to obtain the example trajectories presented in the next subsection.

\subsection{Example Trajectories}
\begin{figure}[h]
     \begin{center}
       \includegraphics[width=10cm]{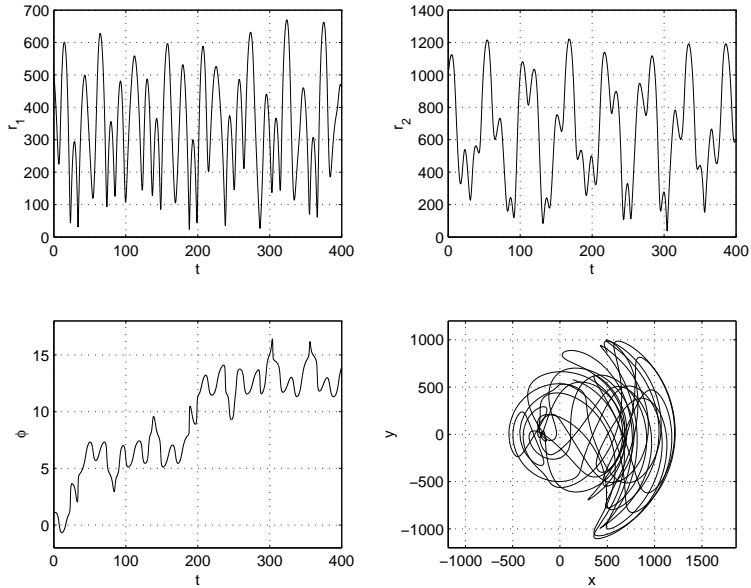}
     \end{center}
     \caption{Time evolution of $r_1$, $r_2$, $\varphi$, and the
trajectory of $m_2$ relative to $m_1$ using RSPRK.}
       \label{fig:rsprk}
\end{figure}

We have computed the reduced trajectory of the double spherical
pendulum using the fourth-order RSPRK algorithm on the Routh
equations. We first consider the evolution of $r_1$,
$r_2$, and $\varphi$, using the RSPRK algorithm on the Routh
equations, as well as the projection of the relative position of
$m_2$ with respect to $m_1$ onto the $xy$ plane as seen in Figure
\ref{fig:rsprk}.

\begin{figure}[h]
\begin{center}
\includegraphics[width=10cm]{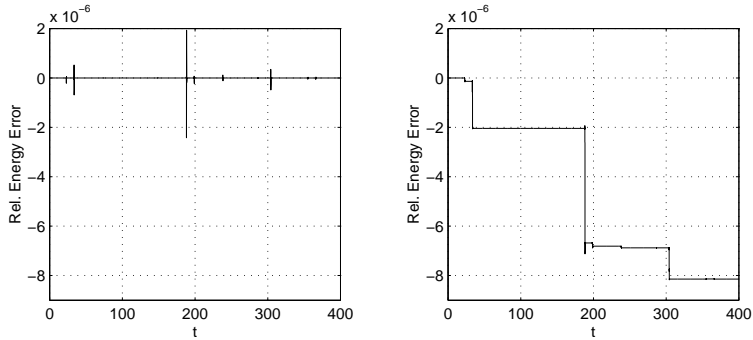}
     \end{center}
     \caption{Relative energy drift $(E-E_0)/E_0$ using RSPRK (left)
compared to the relative energy drift in a non-symplectic RK (right).}
       \label{fig:energy_compare}
\end{figure}

Figure \ref{fig:energy_compare} illustrates that the energy
behavior of the trajectory is very good, as is typical of
variational integrators, and does not exhibit a spurious drift. In
contrast, the non-symplectic fourth-order Runge-Kutta 
applied to the unreduced dynamics has a
systematic drift in the energy, even when using time-steps that are smaller by a factor of four.

\section{Computational Considerations}\label{s-computational}
\medskip 

\noindent {\bf Reduced vs. Unreduced Simulations.}
The reduced dynamics can either be computed directly, by using the discrete Routh or RSPRK equations, or by computing in the unreduced space, and projecting onto the shape space. We discuss the relative merits of these approaches.

Given a configuration space and symmetry group of dimensions $n$ and $m$, respectively, we can either use a simpler algorithm in $2n$ dimensions, or a more geometrically involved algorithm in $2(n-m)$ dimensions. The reduced algorithm involve curvature terms that need to be symbolically precomputed, but these do not affect the sparsity of the system of equations, and the lower dimension results in computational saving that are particularly evident in the repeated, or long-time simulation of problems with a shape space of high codimension.

An example which is of current engineering interest is the dynamics of connected networks of systems
with their own internal symmetries, such as coordinated clusters of satellites modeled as rigid bodies with internal rotors. If the systems to be connected are all identical, the geometric quantities that need to be computed, such as the mechanical connection, have a particularly simple repeated form, and the
small additional upfront effort in implementing the reduced algorithm can result in substantial computational savings due to the lower dimension of the reduced system.

\medskip 

\noindent {\bf Intrinsic vs.~Non-intrinsic methods.} Non-intrinsic numerical schemes, such as SPRK applied to the classical Routh equations, can have undesirable numerical properties due to the need for coordinate dependent local trivializations and the presence of coordinate singularities in these local trivializations, as is the case when using Euler angles for rigid body dynamics.

In the case of non-canonical symplectic forms, frequent computationally expensive coordinate changes are necessary when using standard non-intrinsic schemes, as documented in \cite{WiPeMi1984, Pa1991}, due to the  need to repeatedly apply Darboux' theorem to put the symplectic structure into canonical form. In contrast, intrinsic methods do not depend on a particular choice of coordinate system, and allow for the use of global charts through the use of containing vector spaces with constraints enforced using Lagrange multipliers.

Coordinate singularities can affect the quality of the simulation in subtle ways that may depend on the numerical scheme. In the simulation of the double spherical pendulum, we notice spikes in the energy corresponding to times when $r_1$ or $r_2$ are close to $0$. These errors accumulate in the non-symplectic method, but remain well-behaved in the symplectic method. Alternatively, sharp spikes can
be avoided altogether by evolving the equations on a constraint surface in $\mathbb{R}^3\times\mathbb{R}^3$, as opposed to choosing local coordinates that automatically satisfy the constraints. Here, the increased cost of working in the containing linear space with constraints is offset by not having to transform between charts of $S^2_{l_1}\times S^2_{l_2}$, which can be
significant if the trajectories are chaotic. An extensive discussion of the issue of representations and parametrizations of rotation groups and its implications for computation can be found in \cite{BoTrBo2000}.

Methods which exhibit local conservation properties on each chart, may still exhibit a drift in the conserved quantity across coordinate changes. As discussed in \cite{BeChFa2001}, only methods in which the local representatives of the algorithm commute with coordinate changes exhibit geometric conservation properties that are robust. In particular, intrinsic methods retain their conservation properties through coordinate changes.

Another promising approach is to use the exponential map to update the numerical solution, which is the basis of Lie group integrators, see~\cite{IsMKNoZa2000}. The Lie group approach can yield rigid body integrators that are embedded in the space of $3\times 3$ matrices but automatically evolve on the rotation group, without the use of constraints or reprojection. In \cite{Kr2004, KrEn2004}, analogues of the explicit Newmark and midpoint Lie algorithm are presented that automatically stay on the rotation group, and exhibit good energy and momentum stability. A method based on generating functions and the exponential on Riemannian manifolds is introduced in \cite{LePa1996}. In the variational setting, a Lie group variational integrator for rigid body dynamics is described in \cite{LeLeMc2005}, and higher-order generalizations can be found in \cite{Leok2004b}. 

\section{Conclusions and Future Work}\label{s-conc}

This paper derives the Discrete Routh equations on the discrete shape space
$S\times S$, which are symplectic with respect to a non-canonical
symplectic form, and retains the good energy behavior typically
associated with variational integrators. Furthermore, when the
group action can be expressed as addition, we obtain the Reduced
Symplectic Partitioned Runge-Kutta algorithm on $T^*S$, that can
be considered as a discrete analogue of cotangent bundle
reduction.  By providing an understanding of how the
reduced and unreduced formulations are related at a discrete
level, we enable the user to freely choose whichever formulation
is most appropriate, and provides the most insight to the problem
at hand.

Certainly one of the obvious things to do in the future is to
extend discrete reduction to the case of nonabelian symmetry
groups following the nonabelian version of Routh reduction given
in \cite{JaMa2000, MaRaSc2000}. There are also several
problems, including the averaged $J _2$ problem, in which one can
carry out discrete reduction by stages and relate it to the semidirect product work of \cite{BoSu1999a}. This
is motivated by the fact that the semidirect product reduction
theory of \cite{HoMaRa1998} is a special case of reduction by
stages (at least without the momentum map constraint), as was
shown in \cite{CeHoMaRa1998}.

In further developing discrete
reduction theory, the discrete theory of connections on principal
bundles developed in \cite{LeMaWe2003} is particularly relevant,
as it provides an intrinsic method of representing the reduced
space $(Q\times Q)/G$ as $(S\times S)\oplus\tilde{G}$, where $\tilde{G}=Q\times_G G$ with $G$ acting by conjugation on $G$. Indeed, such a construction can be viewed as a connection on a Lie groupoid, and it is natural to express discrete mechanics on $Q\times Q$ in the language of pair groupoids, as originally proposed in \cite{We1996}. Generalizations of this approach to arbitrary Lie groupoids, as well as a discussion of the role of discrete connections in yielding a discrete analogue of the Lagrange-Poincar\'e equations, can be found in \cite{MaMaMa2005}.

Another component that is needed in this work is a good discrete
version of the calculus of differential forms. Note that in our
work we found, being directed by mechanics, that the right
discrete version of the magnetic $2$-form is the difference of two
connection $1$-forms. It is expected that we could recover such a
magnetic $2$-form by considering the discrete exterior derivative
of a discrete connection form in a finite discretization of
space-time, and taking the continuum limit in the spatial
discretization. Developing a discrete analogue of Stokes' Theorem
would also provide insight into the issue of discrete geometric
phases. Some preliminary work on a discrete theory of exterior
calculus can be found in \cite{DeHiLeMa2003}.

\ack
We gratefully acknowledge helpful comments and suggestions of Alan Weinstein and the referees. SMJ supported in part by ISRO and DRDO through the Nonlinear Studies Group, Indian Institute of Science, Bangalore. ML supported in part by NSF Grant DMS-0504747, and a faculty grant and fellowship from the Rackham Graduate School, University of Michigan. JEM supported in part by NSF ITR Grant
ACI-0204932.


\section*{References}
\begin{thebibliography}{40}



\bibitem{BeChFa2001} Benettin, G., A.~M. Cherubini, and F. Fass\`o, A changing-chart {S}ymplectic algorithm for rigid bodies and other {H}amiltonian systems on manifolds, {\em SIAM J. Sci. Comput.}, \textbf{23}(4), 1189--1203.
 
\bibitem{BoLoSu1998}
Bobenko, A., B.~Lorbeer, and Y.~Suris [1998],
Integrable discretizations of the {E}uler top, {\em Journal of Mathematical
Physics} \textbf{39}, 6668--6683.

\bibitem{BoSu1999}
Bobenko, A. and Y.~Suris [1999], Discrete {L}agrangian reduction,
discrete {E}uler-{P}oincar\'e equations, and semidirect products, {\em Letters
in Mathematical Physics} \textbf{49}, 79--93.


\bibitem{BoSu1999a}
Bobenko, A.~I. and Y.~B. Suris [1999a], Discrete time {L}agrangian mechanics on
    {L}ie groups, with an application to the {L}agrange top, {\em Comm. Math.
    Phys.} \textbf{204}, 147--188.
    
\bibitem{BoTrBo2000}Borri, B., L.~Trainelli, and C.~L.~Bottasso [2000], On representations and parametrizations of motion, {\em Multibody System Dynamics} \textbf{4}, 129--193.


\bibitem{CeHoMaRa1998}
Cendra, H., D.~Holm, J.~Marsden, and T.~Ratiu [1998], Lagrangian Reduction, the
    {E}uler--{P}oincar\'e Equations and Semidirect Products, {\em Amer. Math.
    Soc. Transl.} \textbf{186}, 1--25.

\bibitem{ChMa2003}
Chang, D.~E. and J.~E. Marsden [2003], Geometric derivation of the Delaunay
    variables and geometric phases, {\em Celestial Mechanics and Dynamical
    Astronomy} \textbf{86}, 185--208.


\bibitem{DeHiLeMa2003}
Desbrun, M., A.~Hirani, M.~Leok, and J.~Marsden [2003], Discrete Exterior
    Calculus, arXiv:math.DG/0508341.






\bibitem{HaLuWa2001}
Hairer, E., C.~Lubich, and G.~Wanner [2001], {\em Geometric Numerical
    Integration}.
\newblock Springer, Berlin-Heidelberg-New York.

\bibitem{HaNoWa1993}
Hairer, E., S.~{N\o{}rsett}, and G.~Wanner [1993], {\em Solving Ordinary
    Differential Equations {I} : {N}onstiff problems}, volume~8 of {\em Springer
    Series in Computational Mathematics}.
\newblock Springer-Verlag, second edition.

\bibitem{HaWa1996}
Hairer, E. and G.~Wanner [1996], {\em Solving Ordinary Differential Equations
    {II} : {S}tiff and differential-algebraic problems}, volume~14 of {\em
    Springer Series in Computational Mathematics}.
\newblock Springer-Verlag, second edition.

\bibitem{HoLy2002}
Holm, D.~D and P. Lynch [2002], Stepwise precession of the resonant swinging
spring, {\em SIAM J. Appl. Dynamical Systems}, {\bf 1}, 44--64.

\bibitem{HoMaRa1998}
Holm, D.~D., J.~E. Marsden, and T.~S. Ratiu [1998], The {E}uler--{P}oincar\'{e}
    equations and semidirect products with applications to continuum theories,
    {\em Adv. in Math.} \textbf{137}, 1--81.


\bibitem{IsMKNoZa2000}Iserles, A., H.~Munthe-Kaas, S.~P. N\o rsett, and A.~Zanna [2000], Lie-group methods, {\em Acta Numerica} \textbf{9}, 215--365.


\bibitem{JaMa2000}
Jalnapurkar, S. and J.~Marsden [2000], Reduction of {H}amilton's variational
    principle, {\em Dynamics and Stability of Systems} \textbf{15}, 287--318.


\bibitem{Kr2004} Krysl, P. [2004], Explicit Momentum-conserving Integrator for  Dynamics of Rigid Bodies Approximating the Midpoint Lie Algorithm, {\em Int. J. Numer. Meth. Eng.}, to appear.

\bibitem{KrEn2004}Krysl, P., L.~Endres [2004], Explicit {N}ewmark/{V}erlet algorithm for Time Integration of the Rotational Dynamics of Rigid Bodies, {\em Int. J. Numer. Meth. Eng.}, to appear.


\bibitem{LeLeMc2005}Lee, T.~Y., M.~Leok, N.~H. McClamroch [2005], A Lie Group Variational Integrator for the Attitude Dynamics of a Rigid Body with Applications to the 3D Pendulum, {\em Proc. IEEE Conf. on Control Applications}, 962--967.

\bibitem{LePa1996}Leimkuhler, B. and G.~Patrick [1996], A symplectic integrator for Riemannian manifolds, {\em J. Nonl. Sci.} \textbf{6}, 367--384.

\bibitem{Leok2004}
Leok, M. [2004], {\em Foundations of Computational Geometric Mechanics},
Thesis, California Institute of Technology.

\bibitem{Leok2004b} Leok, M. [2004b], Generalized Galerkin Variational Integrators, arXiv:math.NA/0508360.

\bibitem{LeMaWe2003}
Leok, M., J.~Marsden, and A.~Weinstein [2004], A Discrete Theory of
    Connections on Principal Bundles, arXiv:math.DG/0508338.

\bibitem{LeMaOrWe2003}
Lew, A., J.~E. Marsden, M.~Ortiz, and M.~West [2003], Asynchronous variational
    integrators, {\em Archive for Rat. Mech. An.} \textbf{167}, 85--146.


\bibitem{MaMaMa2005}
Marrero, J.~C., D.~Mart\'{\i}n de Diego, and E.~Mart\'{\i}nez [2005], Discrete Lagrangian and Hamiltonian Mechanics on Lie groupoids, arXiv:math.DG/0506299.

\bibitem{Marsden1992}
Marsden, J. [1992], {\em Lectures on Mechanics}, volume 174 of {\em London
    Mathematical Society Lecture Note Series}.
\newblock Cambridge University Press.

\bibitem{MaMoRa1990}
Marsden, J., R.~Montgomery, and T.~Ratiu [1990], Reduction, Symmetry, and
    Phases in Mechanics, {\em Memoirs of the American Mathematical Society}
    \textbf{88}, 1--110.

\bibitem{MaPeSh1999}
Marsden, J., S.~Pekarsky, and S.~Shkoller [1999], Discrete {E}uler-{P}oincar\'e
and {L}ie-{P}oisson Equations, {\em Nonlinearity} \textbf{12}, 1647--1662.

\bibitem{MaPeSh2000}
Marsden, J., S.~Pekarsky, and S.~Shkoller [2000], Symmetry
    reduction of discrete {L}agrangian mechanics on {L}ie groups,
{\em Journal of Geometry and Physics} \textbf{36}, 140--151.

\bibitem{MaRa1999}
Marsden, J. and T.~Ratiu [1999], {\em Introduction to Mechanics and Symmetry},
    volume~17 of {\em Texts in Applied Mathematics}, Springer-Verlag,
second edition.

\bibitem{MaRaSc2000}
Marsden,
J., T.~Ratiu, andJ.~Scheurle [2000],
Reduction theory and the
{L}agrange-{R}outh equations,
{\em J. Math. Phys.} \textbf{41},
3379--3429.


\bibitem{MaRaWe1984}
Marsden, J.~E., T.~S. Ratiu, and A.~Weinstein [1984], Semi-direct products and
    reduction in mechanics, {\em Trans. Amer. Math. Soc.} 
\textbf{281}, 147--177.

\bibitem{MaSc1993a}
Marsden, J. and J.~Scheurle [1993a], Lagrangian reduction and the
    double spherical pendulum, {\em ZAMP} \textbf{44}, 17--43.

\bibitem{MaSc1993b}
Marsden, J. and J.~Scheurle [1993b], The reduced
    {E}uler--{L}agrange equations, {\em Fields Inst. Commun.} \textbf{1},
    139--164.

\bibitem{MaSc1995}
Marsden, J.~E. and J.~Scheurle [1995], Pattern evocation and geometric phases
    in mechanical systems with symmetry, {\em Dyn. and Stab. of Systems}
    \textbf{10}, 315--338.

\bibitem{MaScWe1996}
Marsden, J.~E., J.~Scheurle, and J.~Wendlandt [1996], Visualization of orbits
    and pattern evocation for the double spherical pendulum.
\newblock In Kirchg{\"{a}}ssner, K., O.~Mahrenholtz, and R.~Mennicken, editors,
    {\em ICIAM 95: Mathematical Research}, volume~87, pages 213--232.
Academie Verlag.

\bibitem{MaWe1974}
Marsden, J.~E. and A.~Weinstein [1974], Reduction of symplectic manifolds with
    symmetry, {\em Rep. Math. Phys.} \textbf{5}, 121--130.

\bibitem{MaWe2001}
Marsden, J. and M.~West [2001], Discrete mechanics and variational integrators.
\newblock In {\em Acta Numerica}, volume~10. Cambridge University Press.


\bibitem{OrRa2004}
Ortega, J.-P. and T.~S. Ratiu [2004], {\em Momentum maps and {H}amiltonian
    reduction}, volume 222 of {\em Progress in Mathematics}.
\newblock Birkh\"auser Boston Inc., Boston, MA.


\bibitem{Pa1991}
Patrick, G. [1991], {\em Two Axially Symmetric Coupled Rigid Bodies: Relative
    Equilibria, Stability Bifurcations, and a Momentum Preserving Symplectic
    Integrator}, PhD thesis, University of California at Berkeley.

\bibitem{PrCo1993}
Prussing, J. and B.~Conway [1993], {\em Orbital Mechanics}.
\newblock Oxford.



\bibitem{SaBlMc2004}
Sanyal, A., A.~J.~Bloch, and H.~N.~McClamroch [2004], {\em Dynamics 
of multibody systems in planar motion in a central
gravitational field}, {\it Dyn. Systems, an Intern. J.} (to appear).

\bibitem{SaShMc2003}
Sanyal, A., J.~Shen, and H.~N.~McClamroch [2003], {\em Variational 
integrators for
    mechanical systems with cyclic generalized coordinates}, (preprint), 2003.



\bibitem{We1996}
Weinstein, A. [1996], Lagrangian mechanics and groupoids. In \textit{Mechanics Days (Waterloo, ON, 1992)}, volume 7 of \textit{Fields Institute Communications}, pages 207-231, American Mathematical Society.

\bibitem{WiPeMi1984}
Wisdom, J., S.~Peale, and F.~Mignard [1984], The Chaotic Rotation of Hyperion,
   {\em Icarus} \textbf{58}, 137--152.

\end{thebibliography}
\end{document}